\documentstyle[12pt]{article}

\title{An Elementary Introduction to the Wiener Process and Stochastic 
 Integrals}
\author{Tam\'as Szabados \\ 
 Technical University of Budapest}
\date{August 23, 1994}

\newcommand{\beq}{\begin{equation}}
\newcommand{\eeq}{\end{equation}}
\newcommand{\beqa}{\begin{eqnarray}}
\newcommand{\eeqa}{\end{eqnarray}}
\newcommand{\beqan}{\begin{eqnarray*}}
\newcommand{\eeqan}{\end{eqnarray*}}

\newcommand{\Pb}{\mbox{\bf P}\left\{ }
\newcommand{\lb}{\left\{ }
\newcommand{\rb}{\right\} }
\newcommand{\E}{\mbox{\bf E}}
\newcommand{\Var}{\mbox{\bf Var}}
\newcommand{\ZZ}{\mbox{\bf Z}}
\newcommand{\RR}{\mbox{\bf R}}

\newtheorem{thm}{Theorem}
\newtheorem{lemma}{Lemma}

\begin{document}

\maketitle

\begin{center} {\em Dedicated to P\'al R\'ev\'esz on the occasion of his 
60th birthday } \end{center}

%...........................

\vspace{10 pt}

%---------------------------------------------------------------------------

\begin{abstract}
An elementary construction of the Wiener process is discussed, based on a 
proper sequence of simple symmetric random walks that uniformly 
converge on bounded intervals, with probability $1$. This method is a 
simplification of F.B. Knight's and P. R\'ev\'esz's. The same sequence is  
applied to give elementary (Lebesgue-type) definitions of It\^o and 
Stratonovich sense stochastic integrals and to prove the basic It\^o formula. 
The resulting approximating sums converge with probability $1$. As a 
by-product, new elementary proofs are given for some properties of the Wiener 
process, like the almost sure non-differentiability of the 
sample-functions. The purpose of using elementary methods almost exclusively 
is twofold: first, to provide an introduction to these topics for a 
wide audience; second, to create an approach well-suited for 
generalization and for attacking otherwise hard problems. 
\end{abstract}

%---------------------------------------------------------------------------

\renewcommand{\thefootnote}{\alph{footnote}}  

%---------------------------------------------------------------------------

\tableofcontents

\section{Introduction}

The {\em Wiener process} is undoubtedly one of the most important stochastic 
processes,
both in the theory and in the applications. Originally it was introduced as 
a mathematical model of {\em Brownian motion}, a random zigzag motion of 
microscopic particles suspended in liquid, discovered by the English 
botanist Brown in 1827. An amazing number of first class scientists 
like Bachelier, Einstein, Smoluchowski, Wiener, and L\'evy, to mention just a 
few, contributed to the theory of Brownian motion. In the course of 
the evolution of probability theory it became clear that the Wiener process 
is a basic tool for many limit theorems and also a natural model
of many phenomena involving randomness, like noise, random fluctuations or 
perturbations. 

The Wiener process is a natural model of Brownian motion. It describes a 
random, but continuous motion of a particle, subjected to the influence of a
large number of chaotically moving molecules of the liquid. Any displacement 
of the particle over an interval of time as a sum of many almost independent 
small influences is normally distributed with expectation zero and 
variance proportional to the length of the time interval. Displacements
over disjoint time intervals are independent. 

The most basic types of {\em stochastic integrals} were introduced by K. It\^o 
and R. L. Straton\-ov\-ich as tools for investigating stochastic differential 
equations, that is, differential equations containing random functions. Not 
surprisingly, the Wiener process is one of the corner stones 
the theory of stochastic integrals and differential equations was built on. 

Stochastic differential equations are applied under similar conditions as 
differential equations in general. The advantage of the stochastic model is 
that it can accommodate noise or other randomly changing input and effects,  
which is a necessity in 
many applications. When solving a stochastic differential equation one has to 
integrate a function with respect to the increments of a stochastic process 
like the Wiener process. In such a case the classical methods of integration 
cannot be applied directly because of the ``strange'' behaviour of the 
increments of the Wiener and similar processes. 

A main purpose of this paper is to provide an elementary introduction to the  
aforementioned topics. The discussion of the 
Wiener process is based on a nice, natural construction of P. R\'ev\'esz 
\cite[Section 6.2]{Revesz 1990}, which is essentially a simplified version of 
F.B. Knight's \cite[Section 1.3]{Knight 1981}. 
We use a proper sequence of simple random walks that converge to the Wiener 
process. Then an elementary definition and discussion of stochastic 
integrals is given, based on \cite{Szabados 1990}, which uses the same 
sequence of random walks. 

The level of the paper is (hopefully) available to any good 
student who has taken a usual calculus sequence and an introductory course 
in probability. Our general reference will be W. Feller's excellent, 
elementary textbook \cite{Feller 1968}. Anything that goes beyond the 
material of that book will be discussed here in detail. I would like to
convince the reader that these important and widely used topics are natural 
and feasible supplements to a strong introductory course in probability; this
way a much wider audience could get acquainted with them. However, I have to 
warn the non-expert reader that ``elementary'' is not a synonym of ``easy'' 
or ``short''.

To encourage the reader it seems worthwhile to emphasize a very useful 
feature of elementary approaches: in many cases, elementary methods are 
easier to generalize or to attack otherwise hard problems. 

%--------------------------------------------------------------------------

\section{Random Walks}

The simplest (and crudest) model of Brownian motion is a 
{\em simple symmetric random walk} in one dimension, hereafter 
{\em random walk} for brevity. 

A particle starts from the origin and steps one unit either to the left or 
to the right with equal probabilities $1/2$, in each unit of time. 
Mathematically, 
we have a sequence $X_1, X_2, \ldots $ of independent and identically 
distributed random variables with 
\[ \Pb X_n=1\rb = \Pb X_n=-1\rb = 1/2 ~(n=1,2,\ldots), \]
and the position of the particle at time $n$ (that is, the random walk) is 
given by the partial sums 
\beq S_0=0,~S_n = X_1 + X_2 + \cdots + X_n ~(n=1,2,\ldots). \label{eq-rw1} 
\eeq
The notation $X(n)$ and $S(n)$ will be used instead of $X_n$ and $S_n$ where 
it seems to be advantageous. 

A bit of terminology: a {\em stochastic process} is a collection $Z(t)~(t 
\in T)$ of random variables defined on a {\em sample space} $\Omega $. 
Usually $T$ is a subset of the real line and $t$ is called ``time''. 
An important concept is that of a {\em sample-function}, that 
is, a randomly selected path of a stochastic process. A sample-function 
of a stochastic process $Z(t)$ can be denoted by $Z(t; \omega)$, where 
$\omega \in \Omega $ is fixed, but the ``time'' $t$ is not.  

To visualize the graph of a sample-function of the random walk one 
can use a broken line connecting the vertices $(n, S_n),~n=1,2,\ldots$ 
(Figure~1). This way the sample-functions are extended from the set of the 
non-negative integers to continuous functions on the interval $[0, \infty)$:
\beq S(t) = S_n + (t-n) X_{n+1}~(n \le  t < n+1;~n=0,1,2,\ldots). 
\label{eq-rw2} \eeq

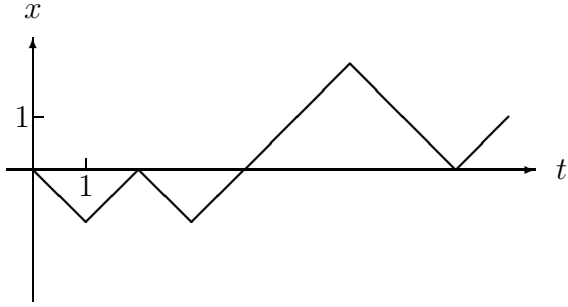
\begin{figure}
\begin{picture}(110,60)(-10,-30)
\thinlines
\put (-5,0){\vector(1,0){100}}
\put (0,-25){\vector(0,1){50}}
\put (10,0){\line(0,1){2}}
\put (0,10){\line(1,0){2}}
\put (100,0){\makebox(0,0){$t$}}
\put (0,30){\makebox(0,0){$x$}}
\put (10,-3){\makebox(0,0){1}}
\put (-2,10){\makebox(0,0){1}}
\thicklines
\put (0,0){\line(1,-1){10}}
\put (10,-10){\line(1,1){10}}
\put (20,0){\line(1,-1){10}}
\put (30,-10){\line(1,1){30}}
\put (60,20){\line(1,-1){20}}
\put (80,0){\line(1,1){10}}
\end{picture}
\caption{The graph of a sample-function of $S(t)$.}  
\end{figure}

It is easy to evaluate the expectation and variance of $S_n$:
\beq \E (S_n) = \sum_{k=1}^n \E (X_k)  = 0, ~~~
 \Var (S_n) = \sum_{k=1}^n \E (X_k^2) = n. \label{eq-rwv} \eeq

The distribution of $S_n$ is a linearly transformed symmetric binomial 
distribution \cite[Section III,2]{Feller 1968}. Each path (broken line) 
of length $n$ has probability $1/2^n$. The number of paths going to the 
point $(n, r)$ from the origin is equal to the number of choosing $(n+r)/2$ 
steps to the right out of $n$ steps. Consequently,
\[ \Pb S_n = r \rb  = {n \choose {(n+r)/2}} {1 \over 2^n}~~(|r| \le n). \] 
The binomial coefficient here is considered to be zero when $n+r$ is not 
divisible by 2. Equivalently, $S_n = 2 B_n - n$, 
where $B_n$ is a symmetric ($p=1/2$) binomial random variable,
$ \Pb B_n = k \rb  = {n \choose k} 2^{-n}$. 

An elementary computation shows that for $n$ large, the binomial distribution 
can be approximated by the normal distribution, see 
\cite[Section VII,2]{Feller 1968}. What is shown there that for  even numbers 
$n=2\nu$ and $r=2k$, if $n \to \infty$ and $| r | < K_n = o(n^{2/3})$, one has
\beq \Pb S_n = r \rb = {n \choose {(n+r)/2}} {1 \over 2^n}
= {2\nu \choose {\nu + k}} {1 \over 2^{2 \nu}} 
\sim {1 \over \sqrt{\pi \nu}} e^{-k^2/\nu}
= 2 h~ \phi ({r h}), \label{eq-bino} \eeq
where $ h = 1/\sqrt n$ and  
$\phi (x) = (1/\sqrt{2 \pi}) e^{-x^2/2}~(-\infty < x < \infty )$,  
the standard normal density function.
Note that for odd numbers $n=2\nu +1$ and $r=2k+1$ (\ref{eq-bino}) can be 
proved similarly as for even numbers.

Here and later we adopt the usual notations
$a_n \sim b_n$  for  $\lim_{n \to \infty} a_n/b_n = 1$  
($a_n$ and $b_n$ are {\em asymptotically equal\/}), and  
$a_n = o(b_n)$  for  $\lim_{n \to \infty} a_n/b_n = 0.$  

Equation~(\ref{eq-bino}) easily implies a special case of the {\em central 
limit theorem} and of the {\em large deviation theorem}, \cite[Sections VII,3 
and 6]{Feller 1968}):

%..........................

\begin{thm} \label{thm-clt1} 
(a) For any real $x$ fixed and $n \to \infty$ we have
\[ \Pb S_n/ \sqrt n~ \le x \rb \to  \Phi (x), \]
where $ \Phi (x) = (1/ \sqrt{2 \pi}) \int_{-\infty}^x e^{-u^2/2}~du $  
  $~(- \infty < x < \infty )$ is the standard normal distribution function. 

(b) If $n \to \infty$ and $x_n \to \infty$ so that $x_n = o(n^{1/6})$, then
\[ \Pb S_n/ \sqrt n~ \ge x_n \rb \sim  1 - \Phi (x_n), \]
\[ \Pb S_n/ \sqrt n~ \le -x_n \rb  \sim \Phi (-x_n) = 1 - \Phi (x_n). 
\Box \]
\end{thm}

%...........................

\vspace{10 pt}

For us the most essential statement of the theorem is that when $x_n$ goes 
to infinity (slower than $n^{1/6}$), then the two sides of (\ref{eq-ldt1}) 
tend to zero equally fast, in fact very fast.
For, to estimate $1 - \Phi(x)$ for $x$ large, one can use the following 
inequality, see \cite[Section VII,1]{Feller 1968}, 
\beq 1 - \Phi (x) < {1 \over {x \sqrt{2 \pi}}} e^{-x^2/2}~~(x > 0). 
\label{eq-asy} \eeq

Thus fixing an $\epsilon > 0$, say $\epsilon = 1/2$, there exists an integer 
$n_0>0$ such that 
\beq \Pb {\left| S_n \over \sqrt n \right|} \ge x_n \rb \le {2 (1+\epsilon ) 
\over x_n \sqrt {2 \pi }} e^{-x_n^2 /2} \le e^{-x_n^2 /2}, 
\label{eq-ldt1} \eeq
for $n \ge n_0$, whenever $x_n \to \infty$ and $x_n = o(n^{1/6})$ as $n \to 
\infty$. It is important to observe that though $S_n$ can take on every 
integer from $-n$ to $n$ with positive probability, the event $\lb |S_n| > 
x_n \sqrt n \rb $ is negligible as $n \to \infty $. 

But what can we do if $n$ does not go to $\infty$, or if the condition 
$x_n=o(n^{1/6})$ does not hold? Then a simple, but still powerful tool, 
{\em Chebyshev's inequality} can be used. A standard form of Chebyshev's 
inequality \cite[Section IX,6]{Feller 1968} is
\[ \Pb | X - \E (X) | \ge t \rb \le {\Var (X) \over t^2} , \]
for any $t > 0$, supposing $\Var (X)$ is finite. An other form that can be 
proved similarly is 
\beq \Pb | X | \ge t \rb \le {\E (|X|) \over t} , \label{eq-che1} \eeq
for any $t > 0$ if $\E (X)$ is finite. If the {\em $k$th moment} of $X$, 
$\E (X^k)$ is finite $(k > 0)$, then one can apply (\ref{eq-che1}) to 
$|X|^k$ getting 
\[ \Pb |X| \ge t \rb = \Pb |X|^k \ge t^k \rb \le {\E (|X|^k) \over t^k}, \]
for any $t > 0$. 

One can even get an upper bound going to $0$ exponentially 
fast as $t \to \infty$ if $\E (e^{uX})$, the {\em moment generating 
function of $X$}, is finite for some $u_0 > 0$. For then, by (\ref{eq-che1}),    
\beq \Pb X \ge t \rb = \Pb u_0 X \ge u_0 t \rb = \Pb e^{u_0 X} \ge e^{u_0 t} 
\rb \le e^{-u_0 t} \E (e^{u_0 X}), 
\label{eq-chee1} \eeq 
for any $t > 0$. 

Analogously, if $\E (e^{-u_0 X})$ is finite for some $u_0 > 0$, then 
\beq \Pb X \le -t \rb = \Pb -u_0 X \ge u_0 t \rb = \Pb e^{-u_0 X} \ge 
e^{u_0 t} \rb \le e^{-u_0 t} \E (e^{-u_0 X}), 
\label{eq-chee2} \eeq
for any $t > 0$. Combining (\ref{eq-chee1}) and (\ref{eq-chee2}), 
one gets
\beq \Pb |X| \ge t \rb = \Pb X \ge t \rb + \Pb X \le -t \rb \le e^{-u_0 t} 
\left( \E (e^{u_0 X}) + \E (e^{-u_0 X}) \right), 
\label{eq-chee} \eeq
for any $t > 0$ if the moment generating function is finite both at $u_0$ 
and at $-u_0$. 

Now, it is easy to find the moment generating function of one step of the 
random walk: 
\[ \E (e^{u X_k}) = e^u (1/2) + e^{-u} (1/2) = \cosh u . \]
Hence, using the independence of the steps, one obtains the moment generating 
function of the random walk $S_n$ as 
\beq \E (e^{u S_n}) = \E ( e^{u \sum _{k=1}^n X_k}) = \E (\prod _{k=1}^n 
e^{u X_k}) = (\cosh u)^n ~~(-\infty < u < \infty ,~n \ge 0).
\label{eq-mgfsn} \eeq

Since $\cosh u$ is an even function and $\cosh 1 < 2$, (\ref{eq-chee}) implies 
that
\beq \Pb |S_n| \ge t \rb \le 2 \cdot 2^n e^{-t} ~~(t>0, ~n \ge 0). 
\label{eq-chesn} \eeq

%-------------------------------------------------------------------------

\section{Waiting Times} 

In the sequel we need the distribution of the random time $\tau$ when a 
random walk first hits either the point $x=2$ or $-2$:
\beq \tau = \tau _1 = \min \lb n : | S_n | = 2 \rb . \label{eq-tau} \eeq
To find the probability distribution of $\tau$, imagine the random walk as a 
sequence of pairs of steps. These (independent) pairs can be classified 
either as a ``return'': $(1, -1)$ or $(-1, 1)$, or as a ``change of magnitude 
2'': $(1, 1)$ or $(-1, -1)$. Both cases have the same probability $1/2$. 

Clearly, it has zero probability that $\tau$ is equal to an odd number. The 
event $\lb \tau = 2j \rb$ occurs exactly when $j-1$ ``returns'' are followed 
by a ``change of magnitude 2''. Because of the independence of the pairs of 
steps, $\Pb \tau = 2j \rb = 1/2^j $.
It means that $\tau = 2 Y$, where $Y$ has geometric distribution with 
parameter $p=1/2$,
\beq \Pb \tau = 2j \rb = \Pb Y = j \rb = 1/2^j~~(j \ge 1). 
\label{eq-geom} \eeq
Hence,
\beq \E (\tau) = 2 \E (Y) = 2(1/p) = 4,~~~ \Var (\tau) = 2^2 
\Var (Y) = 2^2 (1-p)/p^2 = 8. \label{eq-etau} \eeq

An important consequence is that with probability $1$, a random walk sooner 
or later hits $2$ or $-2$:
\[ \Pb \tau < \infty \rb = \sum_{j=1}^{\infty } (1/2^j) = 1. \]

It is also quite obvious that
\beq \Pb S(\tau) = 2 \rb = \Pb S(\tau) = -2 \rb = 1/2. \label{eq-ru1} \eeq
This follows from the symmetry of the random walk. If we reflect 
$S(t)$ to the time axis, the resulting process $S^*(t)$ is also 
a random walk. Its corresponding $\tau ^*$ is equal to $\tau $, and the event 
$\lb S^*(\tau ) = 2 \rb $ is the same as $\lb S(\tau ) = -2 \rb $. Since 
$S^*(t)$ is just the same sort of random walk as $S(t)$, we have 
$\Pb S^*(\tau )=2 \rb = \Pb S(\tau )=2 \rb $ as well. 

Another way to show (\ref{eq-ru1}) is to use the fact that the waiting time 
$\tau $ has countable many possible values and for any specific value we have 
symmetry:  
\beqan
\Pb S(\tau ) = 2 \rb 
& = & \sum _{j=1}^{\infty} \Pb S(2j) = 2 \mid \tau = 2j \rb \Pb \tau = 2j \rb 
 \\
& = & \sum _{j=1}^{\infty} \Pb A_{2j-2}, X_{2j} = X_{2j-1} = 1 \mid 
 A_{2j-2}, X_{2j} = X_{2j-1} \rb \Pb \tau = 2j \rb \\
& = & (1/2) \sum _{j=1}^{\infty} \Pb \tau = 2j \rb = 1/2, 
\eeqan
where $A_{2j-2}$ denotes the event that each of the first $j-1$ pairs is a 
``return'', i.e. $A_{2j-2}=\lb X_2 = -X_1, \ldots , X_{2j-2} = -X_{2j-3} \rb$, 
$A_0=\emptyset $. 

We mention that (\ref{eq-ru1}) illustrates a consequence of the so-called 
optional sampling theorem too: $\E S(\tau )) = 2 \Pb S(\tau )=2 \rb + (-2) 
\Pb S(\tau )=-2 \rb = 0$, which is the same as the expectation of $S(t)$. 

We also need the probability of the event that a random walk starting from 
the point $x=1$ hits $x=2$ before hitting $x=-2$. This is equal to the 
conditional probability $\Pb S(\tau ) = 2 \mid X_1 = 1 \rb $.  
If $X_1 = 1$, then $X_2 = 1$ with probability $1/2$, and then $\tau =2$ and 
$S(\tau )=2$ as well:   
$ \Pb S(\tau ) = 2, \tau = 2 \mid X_1 = 1 \rb = 1/2. $ 

On the  other hand, if $X_1 = 1$, then $\tau > 2$ if and only if $X_2=-1$, 
with probability $1/2$. that is, at the second step the walk returns the 
origin and starts ``from scratch''. Then by (\ref{eq-ru1}), it has probability 
$1/2$ that the random walk hits $2$ sooner than $-2$:  
$ \Pb S(\tau ) = 2, \tau > 2 \mid X_1 = 1 \rb = 1/4$. 
Therefore
\beqa 
\Pb S(\tau ) = 2 \mid X_1 = 1 \rb 
& = & \Pb S(\tau ) = 2, \tau = 2 \mid X_1 = 1 \rb + \Pb S(\tau ) = 2, 
 \tau > 2 \mid X_1 = 1 \rb \nonumber \\
& = & (1/2) + (1/4) = 3/4. \label{eq-ru2} 
\eeqa
It also follows then that 
\beq \Pb S(\tau ) = -2 \mid X_1 = 1 \rb = 1 - (3/4) = 1/4. 
\label{eq-ru3} \eeq
(\ref{eq-ru1}), (\ref{eq-ru2}), and (\ref{eq-ru3}) are special cases of 
ruin probabilities \cite[Section XIV,2]{Feller 1968}. For example, it can   
be shown that the probability that a random walk hits the level $a > 0$ 
before hitting the level $-b < 0$ is $b/(a+b)$. 

Extending definition (\ref{eq-tau}) of $\tau $, for $k = 1, 2,\ldots $ we 
recursively define
\[ \tau _{k+1} = \min \lb n : n > 0, | S(T_k + n) - S(T_k) |  = 2 \rb , \]
where
\beq T_k = T(k) = \tau_1 + \tau_2 + \cdots + \tau_k. \label{eq-nebi} \eeq

Then each $\tau_k$ has the same distribution as $\tau = \tau_1 $. For, 
\beqan
\lefteqn
{ \Pb \tau_{k+1} = 2j \mid T_k = 2m \rb } \\
& = & \Pb ~\min \lb n : n > 0,  | S(2m+n) - S(2m) | = 2 \rb~ = 2j 
 \mid T_k = 2m \rb \\
& = & \Pb ~\min \lb n : n > 0,  | S(n)| = 2 \rb~ = 2j \rb 
 = \Pb \tau_1 = 2j \rb = 1/2^j, 
\eeqan
where $k \ge 1$, $j \ge 1$, and $m \ge 1$ are arbitrary.
The second equality above follows from two facts. First, each increment 
$S(2m+n) - S(2m)$ is independent of the event $\lb T_k = 2m \rb $, because 
the increment depends only on the random variables $X_i~ (2m+1 \le i \le 
2m+n)$, while the event $\lb T_k = 2m \rb $ is determined exclusively by the 
random variables $X_i~ (1 \le i \le 2m)$, the corresponding ``past''. Second, 
each increment $S(2m+n) - S(2m)$ has the same distribution as $S(n)$, since 
both of them is a sum of $n$ independent $X_i$.
Hence, $\tau_{k+1}$ is independent of $T_k$ (and also of any $\tau_i, i \le 
k$), so indeed,
$ \Pb \tau_{k+1} = 2j \rb  = 1/2^j$ $~(j \ge 1)$. 

%...........................

\vspace{10 pt}

We also need the distribution of the random time $T_k$ required by $k$ 
changes of magnitude 2 along the random walk.
In other words, $S(t)$ hits even integers (different from the previous one) 
exclusively at the time instants $T_1, T_2, \ldots $.
To find the probability distribution of $T_k$, imagine the random walk again 
as a sequence of independent pairs of steps, ``returns'' and ``changes of 
magnitude 2'', both types having probability $1/2$.
The number of cases the event $\lb T_k = 2j \rb ~(j \ge k)$ can occur 
is equal to the number of choices of $k-1$ pairs out of $j-1$ where 
a change of magnitude 2 occurs, before the last pair, which is necessarily a 
change of magnitude 2. Therefore 
\beq \Pb T_k = 2j \rb = {j-1 \choose k-1} {1 \over 2^j}~~~(j \ge k \ge 1). 
\label{eq-nbd} \eeq
It means that $T_k = 2 N_k$, where $N_k$ has a negative binomial 
distribution 
with $p=1/2$, \cite[Section VI,8]{Feller 1968}. 

All this also follows from the 
fact that $N_k=T_k/2$ is the sum of $k$ independent, geometrically 
distributed random variables with parameter $p=1/2$, see (\ref{eq-geom}) and 
(\ref{eq-nebi}): $N_k = Y_1 + Y_2 + \cdots + Y_k~(Y_j = \tau_j/2).$
Then $T_k$ is finite valued with probability $1$ and 
the expectation and variance of $T_k$ easily follows from (\ref{eq-etau}) 
and (\ref{eq-nebi}): 
\beq \E (T_k) = k \E (\tau) = 4k, ~~~\Var (T_k) = k \Var (\tau) = 8k. 
\label{eq-nbev} \eeq 

It is worth mentioning that $T_k$ is a {\em stopping time} for each $k \ge 
1$. By definition, it means that any event of the form $\lb T_k \le j \rb $ 
depends exclusively on the corresponding ``past'' $S(t)~(t \le j)$. In other 
words, $S_1, \ldots , S_j$ determine whether $\lb T_k \le j \rb $ occurs or 
not.  

Fortunately, the central limit and the large deviation theorems (see Theorem 
\ref{thm-clt1}) can be proved for negative binomial distributions in the same 
fashion as for binomial distributions. 

%................................

\begin{thm} \label{thm-clt2}
(a) For any real $x$ fixed and $k \to \infty$ we have  
\[ \Pb {T_k - 4k \over \sqrt{8k}} \le x \rb \to \Phi (x). \]

(b) If $k \to \infty$ and $x_k \to \infty$ so that $x_k = o(k^{1/6})$, then
\[ \Pb {T_k - 4k \over \sqrt{8k}} \ge x_k \rb \sim 1 - \Phi (x_k), \]
\[ \Pb {T_k - 4k \over \sqrt{8k}} \le -x_k \rb \sim \Phi(-x_k) = 1 - 
\Phi (x_k). \]
\end{thm}

{\sc Proof.}
The normal approximation (\ref{eq-bino}) is applicable to negative binomial                               
distributions too: if $r=2j$ and $k \to \infty$, then  
\beqa
\Pb T_k = r \rb & = & {j-1 \choose k-1} {1 \over 2^j}
 = {1 \over 2} {j-1 \choose {(j-1)+(2k-j-1) \over 2}} {1 \over 2^{j-1}} 
 \nonumber \\
& \sim & {1 \over 2} {1 \over \sqrt{\pi (j-1)/2}} \exp \left( -{({2k-j-1 
 \over 2})^2 \over (j-1)/2} \right) \nonumber \\
& = & {1 \over \sqrt{\pi (r-2)}} \exp \left( -{(r-4k+2)^2 \over 4r-8} \right), 
\label{eq-nbn} \eeqa
supposing $ |2k-j-1| = o((j-1)^{2/3})$, or equivalently, 
\beq |r-4k| = o(k^{2/3}). \label{eq-nbo} \eeq

A routine computation shows that (\ref{eq-nbn}) is asymptotically equal to 
\[ \sim {1 \over \sqrt{4k\pi}} \exp \left( -{(r-4k)^2 \over 16k} \right) , \]
when $k \to \infty $ and (\ref{eq-nbo}) holds. Therefore we get an analogue of 
(\ref{eq-bino}): if $k \to \infty$ and $r$ is any even number such that 
$|r - 4k| < K_k = o(k^{2/3}),$
\beq \Pb T_k = r \rb \sim 2h~\phi((r-4k)h),~ h=1/\sqrt{8k}, \label{eq-nbno} 
\eeq
where $\phi$ denotes the standard normal density function. 

Then in the same way as the statements of Theorem \ref{thm-clt1} are 
obtained from (\ref{eq-bino}) in \cite[Sections VII,3 and 6]{Feller 1968} one 
can get the present theorem from (\ref{eq-nbno}). Here we recall only the 
basic step of the argument: 
\beqan
\Pb x_1 \le {T_k - 4k \over \sqrt {8k}} \le x_2  \rb 
& \sim & \sum _{\lb r:~x_1 \le (r-4k)h \le x_2,~r~\mbox{\scriptsize is even} 
\rb } 2h~\phi((r-4k)h) \\
& \to & \int_{x_1}^{x_2} \phi(t)~dt  = \Phi (x_2) - \Phi (x_1), 
\eeqan
for any $x_1, x_2$, when $k \to \infty$ and so $h \to 0$. The simple 
meaning of this is that Riemann sums converge to the corresponding 
integral. $\Box$

%...........................

\vspace{10 pt}

In the same fashion as the large deviation inequality (\ref{eq-ldt1}) 
was obtained for $S_n$, Theorem \ref{thm-clt2}(b) and (\ref{eq-asy}) imply 
a large deviation type inequality for $T_k$:  
\beq \Pb {\left| T_k - 4k \over \sqrt{8k} \right| } \ge x_k \rb
\le e^{-x_k^2 /2}, \label{eq-ldt2} \eeq
for $k \ge k_0$, supposing $x_k \to \infty$ and $x_k = o(k^{1/6})$ as $k \to 
\infty$. 

Like in case of $S_n$, with $T_k$ too we need a substitute for the large 
deviation inequality if the assumptions $k \to \infty$ or $x_k = 
o(k^{1/6})$ do not hold. The moment generating function of $\tau _n$ is 
simple:
\beq \E (e^{u \tau _n}) = \sum _{j=1}^{\infty} e^{u2j} {1 \over 2^j} 
={ e^{2u}/2 \over 1- (e^{2u}/2)} = {1 \over 2 e^{-2u} -1}. \label{eq-mgftau} 
\eeq
This function is finite if $u < \log \sqrt 2$. Here and afterwards $\log$ 
denotes logarithm with base $e$. 

Now the moment generating function of $T_k$ follows from the independence of 
the $\tau _n$'s as 
\beq \E (e^{u T_k}) = \E (e^{u \sum _{n=1}^k \tau _n}) = 
\E (\prod _{n=1}^k e^{u \tau _n}) 
= (2 e^{-2u} - 1)^{-k} ~~(u < \log \sqrt 2 , ~ k \ge 0).
\label{eq-mgftk} \eeq

We also need the moment generating function of the centered and ``normalized'' 
random variable $(T_k -  4k)/\sqrt 8$, whose expectation is $0$ and variance 
is $k$: 
\beq 
\E (e^{u (T_k - 4k)/\sqrt 8}) = e^{-4ku/\sqrt 8} ~\E (e^{T_k u/\sqrt 8}) 
=  \left( 2 e^{u/\sqrt 2} - e^{u \sqrt 2} \right) ^{-k}, 
 \label{eq-mgftkc} \eeq
for $u < \sqrt 2 \log 2$ and $k \ge 0$.
Since (\ref{eq-mgftkc}) is less than $2^k$ for $u= \pm 1/2$, the exponential 
Chebyshev's inequality (\ref{eq-chee}) implies that 
\beq \Pb |T_k - 4k|/\sqrt 8 \ge t \rb \le 2 \cdot 2^k e^{-t/2} ~~(t>0,
 ~ k\ge 0). \label{eq-chetk} \eeq

%--------------------------------------------------------------------------

\section{From Random Walks to the Wiener Process: \protect\newline 
 ``Twist and Shrink''}
\label{se-rwtw}

Our construction of the Wiener process is based on P. R\'ev\'esz's one, \cite
[Section 6.2]{Revesz 1990}, which in turn is a simpler version of F.B. 
Knight's \cite[Section 1.3]{Knight 1981}. The advantage of this method over 
the several known ones is that it is very natural and elementary. 

We will define a sequence of approximations to the Wiener process, each of 
which is a ``twisted and shrunk'' random walk, a refinement of the previous 
one. It will be shown that this sequence converges to a process 
having the properties characterizing the Wiener process. 

Imagine that we observe a particle undergoing Brownian motion. In the first 
experiment we observe the particle exclusively when it hits points with 
integer coordinates $j\in \ZZ $. Suppose that it happens exactly at the 
consecutive time instants 
$1,2,\ldots$. To model the graph of the particle between the vertices 
so obtained the simplest idea is to join them by straight line segments like 
in Figure~1. Therefore the first approximation is  
\[ B_0(t) = S_0(t) = S(t), \]
where $t \ge 0 $ real and $S(t)$ is a random walk defined by (\ref{eq-rw1}) 
and (\ref{eq-rw2}). 

Suppose that in the second experiment we observe the particle when it hits 
points with coordinates $j/2~(j\in \ZZ )$, in the third experiment when it 
hits points with coordinates $j/2^2~(j \in \ZZ )$, etc. To model the second 
experiment one idea is to take a second random walk $S_1(t)$, independent 
of the first one, and shrink it. 

Then the first problem that arises is the relationship 
between the time and space scales: if one wants to compress the length of a 
step into half, how much one has to compress the time needed for one step to 
preserve the 
essential properties of a random walk. Here we recall that by (\ref{eq-rwv}),  
the square root of the average squared distance of the random walk from the 
origin after time $n$ is $\sqrt n$. So shrinking the random walk so that 
there are $n$ steps in one time unit, each step should have a length 
$1/\sqrt n$. This way after one time unit the square root of the average 
squared 
distance of the walk from the origin will be around one spatial unit, like 
in the case of the original random walk. It means that compressing the 
length of one step into $1/2$ (or in general: $1/2^m, m=1,2,\ldots $) 
one has to compress the time needed for one step into $1/2^2$ (in general:  
$1/2^{2m}$).

The second problem is that sample-functions of $B_0(t)$ and of a shrunk version 
of an independent $S_1(t)$ have nothing to do with each other, the second is 
not being a refinement of the first in general. For example, 
if $B_0(1)=1$, then it is equally likely that the first integer the shrunk 
version of $S_1(t)$ hits is $+1$ or $-1$. 

Hence before shrinking we want to 
modify $S_1(t)$ so that it hits even integers $2j~(j \in \ZZ )$ 
(counting the next one only if it is different from the previous one) in 
exactly the same order as $S_0(t)$ hits the corresponding integers
$j \in \ZZ $. For example, if $S_0(1)=1$ and $S_0(2)=2$, then the first
even integer $S_1(t)$ hits should be $2$ and the next one (different from
2) should be $4$. Thus if $S_1(t)$ hits the first even integer at 
time $T_1(1)$ and $S_1(T_1(1))$
happens to be $-2$, we will reflect every step $X_1(k)$ of $S_1(t)$ for
$0 < k \le T_1(1)$. This way we get a modified random walk $\tilde S_1(t)$
up to time $T_1(1)$ so that $\tilde S_1(T_1(1)) = 2$. Then we continue
similarly up to time $T_1(2)$: if the (already modified) walk hit
$0$ at time $T_1(2)$ (instead of $4$), then we would reflect the steps
$X_1(k)$ for $T_1(1) < k \le T_1(2)$.
This modification process, which we will call ``twisting'', ensures that the
next approximation will always be a refinement of the previous one.

\vspace{10 pt}

Now let us see the construction in detail. 
It begins with a sequence of independent random walks $S_0(t), S_1(t),
\ldots $. That is, for each $ m \ge 0 $,
\beq S_m(0) = 0,~S_m(n) = X_m(1) + X_m(2) + \cdots + X_m(n)~(n \ge 1),
\label{eq-smn} \eeq
where $ X_m(k)~(m \ge 0, k \ge 1) $ is a double array of independent,
identically distributed random variables such that
\beq \Pb X_m(k) = 1 \rb = \Pb X_m(k) = -1 \rb = 1/2. \label{eq-xmi} \eeq

First we possibly modify $S_1(t), S_2(t), \ldots $ one-by-one, 
using the ``twist'' method to obtain a sequence of {\em not} 
independent random walks $ \tilde S_1(t), \tilde S_2(t), \ldots $, each of 
which is a refinement of the former one.
Second, by shrinking we get a sequence $B_1(t), B_2(t), \ldots $ 
approximating the Wiener process.

In accordance with the notation in (\ref{eq-nebi}), for $m \ge 1$, $S_m$
hits even integers (different from the previous one)
exclusively at the random time instants  
\[ T_m(0) = 0,~T_m(k) = \tau _m(1) + \tau _m(2) + \cdots + \tau _m(k)~
(k \ge 1). \]
Each random variable $T_m(k)$ has the same distribution as $T(k)=T_k$ above, 
see (\ref{eq-nbd}) and (\ref{eq-nbev}). That is, $T_m(k)$ is the double of a 
negative binomial random variable, with expectation $4k$ and 
variance $8k$.

Now we define a suitable sequence of ``twisted'' random walks $\tilde 
S_m(t)$~ $(m \ge 1)$ recursively, using $\tilde S_{m-1}(t)$, starting with
\[ \tilde S_0(t) = S_0(t)~(t \ge 0). \] 
First we set
\[ \tilde S_m(0) = 0. \]
Then for $k = 0, 1, \ldots$ successively and for every $n$ such that $T_m(k) 
< n \le T_m(k+1)$, we take (Figures~2-4).
\beq \tilde X_m(n) = \left\{ \begin{array}{rl}
X_m(n) & \mbox{ if $~~S_m(T_m(k+1)) - S_m(T_m(k)) 
 = 2 \tilde X_{m-1}(k+1)$;} \\
-X_m(n) & \mbox{ otherwise.} 
\end{array} \right. \label{eq-txd} \eeq
and 
\beq  \tilde S_m(n) = \tilde S_m(n-1) + \tilde X_m(n).
\label{eq-tsmn} \eeq

Observe that the stopping times $\tilde T_m(k)$ corresponding to 
$\tilde S_m(t)$ are the same as the original ones $T_m(k)$ 
~$(m \ge 0, k \ge 0)$. 

%----------------------------------------------------------------

\begin{figure}
\begin{picture}(50,40)(-10,-10)
\thinlines
\put (-5,0){\vector(1,0){40}}
\put (0,-5){\vector(0,1){30}}
\put (10,0){\line(0,1){2}}
\put (20,0){\line(0,1){2}}
\put (30,0){\line(0,1){2}}
\put (0,10){\line(1,0){2}}
\put (0,20){\line(1,0){2}}
\thicklines
\put (40,0){\makebox(0,0){$t$}}
\put (0,30){\makebox(0,0){$x$}}
\put (10,-3){\makebox(0,0){1}}
\put (20,-3){\makebox(0,0){2}}
\put (30,-3){\makebox(0,0){3}}
\put (-2,10){\makebox(0,0){1}}
\put (-2,20){\makebox(0,0){2}}
\put (0,0){\line(1,1){20}}
\put (20,20){\line(1,-1){10}}
\put (0,0){\circle*{3}}
\put (10,10){\circle*{3}}
\put (20,20){\circle*{3}}
\put (30,10){\circle*{3}}
\end{picture}
\caption{$B_0(t;\omega ) = S_0(t;\omega )$.}
\end{figure}
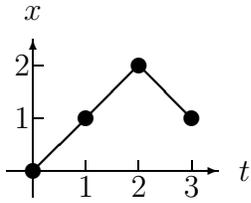

\begin{figure}
\begin{picture}(175,55)(-5,-25)
\thinlines
\put (-5,0){\vector(1,0){170}}
\put (0,-25){\vector(0,1){50}}
\put (10,0){\line(0,1){2}}
\put (0,20){\line(1,0){2}}
\put (0,-20){\line(1,0){2}}
\thicklines
\put (170,0){\makebox(0,0){$t$}}
\put (0,30){\makebox(0,0){$x$}}
\put (10,-3){\makebox(0,0){1}}
\put (-2,20){\makebox(0,0){2}}
\put (-3,-20){\makebox(0,0){-2}}
\put (0,0){\line(1,-1){10}}
\put (10,-10){\line(1,1){10}}
\put (20,0){\line(1,-1){10}}
\put (30,-10){\line(1,1){20}}
\put (50,10){\line(1,-1){30}}
\put (80,-20){\line(1,1){20}}
\put (100,0){\line(1,-1){10}}
\put (110,-10){\line(1,1){20}}
\put (130,10){\line(1,-1){10}}
\put (140,0){\line(1,1){20}}
\put (0,0){\circle*{3}}
\put (80,-20){\circle*{3}}
\put (100,0){\circle*{3}}
\put (160,20){\circle*{3}}
\end{picture}
\caption{$S_1(t;\omega )$.}
\end{figure}
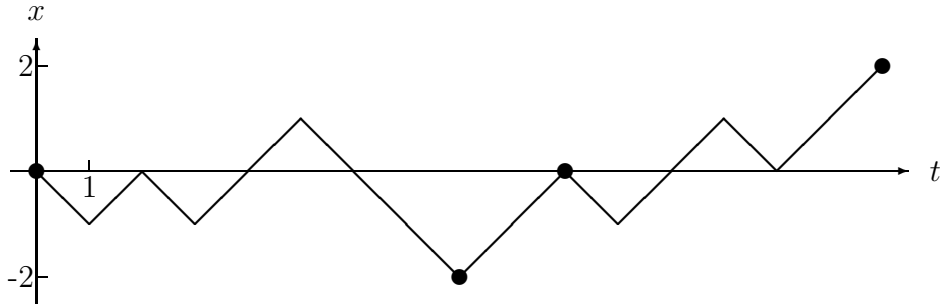

\begin{figure}
\begin{picture}(175,80)(-5,-10)
\thinlines
\put (-5,0){\vector(1,0){170}}
\put (0,-10){\vector(0,1){75}}
\put (10,0){\line(0,1){2}}
\put (80,0){\line(0,1){2}}
\put (100,0){\line(0,1){2}}
\put (160,0){\line(0,1){2}}
\put (0,20){\line(1,0){2}}
\put (0,40){\line(1,0){2}}
\put (0,60){\line(1,0){2}}
\thicklines
\put (170,0){\makebox(0,0){$t$}}
\put (0,70){\makebox(0,0){$x$}}
\put (10,-3){\makebox(0,0){1}}
\put (80,-3){\makebox(0,0){$T_1(1)$}}
\put (100,-3){\makebox(0,0){$T_1(2)$}}
\put (160,-3){\makebox(0,0){$T_1(3)$}}
\put (-2,20){\makebox(0,0){2}}
\put (-2,40){\makebox(0,0){4}}
\put (-2,60){\makebox(0,0){6}}
\put (0,0){\line(1,1){10}}
\put (10,10){\line(1,-1){10}}
\put (20,0){\line(1,1){10}}
\put (30,10){\line(1,-1){20}}
\put (50,-10){\line(1,1){60}}
\put (110,50){\line(1,-1){20}}
\put (130,30){\line(1,1){10}}
\put (140,40){\line(1,-1){20}}
\put (0,0){\circle*{3}}
\put (80,20){\circle*{3}}
\put (100,40){\circle*{3}}
\put (160,20){\circle*{3}}
\end{picture}
\caption{$\tilde S_1(t;\omega )$.}
\end{figure}
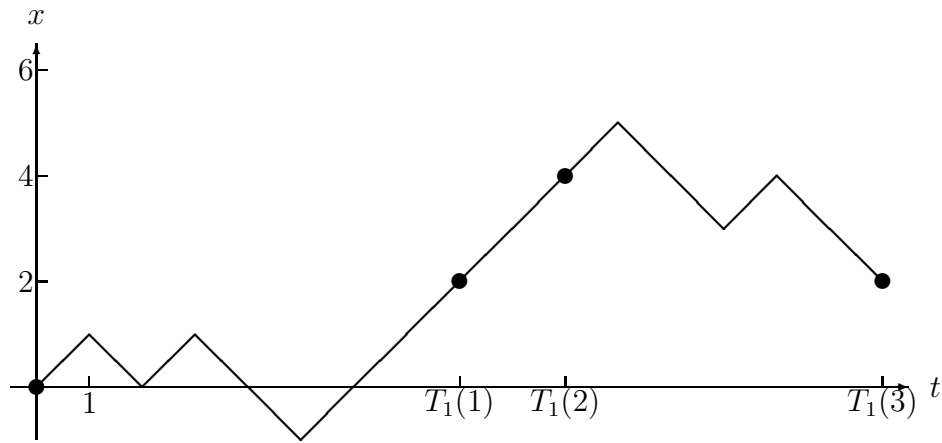

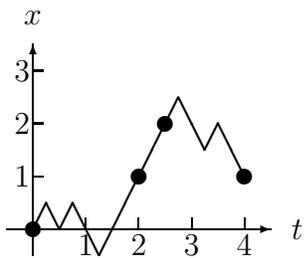
\begin{figure}
\begin{picture}(60,50)(-10,-10)
\thinlines
\put (-5,0){\vector(1,0){50}}
\put (0,-5){\vector(0,1){40}}
\put (10,0){\line(0,1){2}}
\put (20,0){\line(0,1){2}}
\put (30,0){\line(0,1){2}}
\put (40,0){\line(0,1){2}}
\put (0,10){\line(1,0){2}}
\put (0,20){\line(1,0){2}}
\put (0,30){\line(1,0){2}}
\thicklines
\put (50,0){\makebox(0,0){$t$}}
\put (0,40){\makebox(0,0){$x$}}
\put (10,-3){\makebox(0,0){1}}
\put (20,-3){\makebox(0,0){2}}
\put (30,-3){\makebox(0,0){3}}
\put (40,-3){\makebox(0,0){4}}
\put (-2,10){\makebox(0,0){1}}
\put (-2,20){\makebox(0,0){2}}
\put (-2,30){\makebox(0,0){3}}
\put (0,0){\line(1,2){2.5}}
\put (2.5,5){\line(1,-2){2.5}}
\put (5,0){\line(1,2){2.5}}
\put (7.5,5){\line(1,-2){5}}
\put (12.5,-5){\line(1,2){15}}
\put (27.5,25){\line(1,-2){5}}
\put (32.5,15){\line(1,2){2.5}}
\put (35,20){\line(1,-2){5}}
\put (0,0){\circle*{3}}
\put (20,10){\circle*{3}}
\put (25,20){\circle*{3}}
\put (40,10){\circle*{3}}
\end{picture}
\caption{$B_1(t;\omega )$.}
\end{figure}

%-----------------------------------------------------------------

%............................

\begin{lemma} \label{l-rw}
For each $m \ge 0$, $\tilde S_m(t)$ ($t \ge 0$) is a random walk, that is, 
$\tilde X_m(1), \tilde X_m(2), \ldots $ is a sequence of independent, 
identically distributed random variables such that
\beq \Pb \tilde X_m(n) = 1 \rb = \Pb \tilde X_m(n) = -1 \rb = 1/2~~(n \ge 1).
\label{eq-txmi} \eeq
\end{lemma}

{\sc Proof.} 
We proceed by induction over $m \ge 0$. For $m=0$, $\tilde S_0(t) = S_0(t)$,  
a random walk by definition. So assume that $\tilde S_{m-1}(t)$ is a random 
walk, where $m \ge 1$, and see if it implies that $\tilde S_m(t)$ is a random 
walk too.

It is enough to show that for any $n \ge 1$ and any $\epsilon _j = \pm 1 ~
(j = 1, \ldots , n)$ we have
\beq \Pb \tilde X_m(1) = \epsilon _1, \ldots , \tilde X_m(n-1) = 
 \epsilon _{n-1}, \tilde X_m(n) = \epsilon _n \rb = 1/2^n. 
\label{eq-iid} \eeq

Set the events $ A_{m,r} = \lb \tilde X_m(j) = \epsilon _j, ~1 \le j \le r 
\rb $ for $1 \le r \le n$ ($A_{m,0}$ is the sure event by definition) and the 
random variables  
$\Delta S_{m,k}^*=S_m(T_m(k + 1)) - S_m(T_m(k))$ for $k \ge 0$. The event 
$A_{m,n-1}$ determines the greatest integer $k \ge 0$ such that $T_m(k) \le 
n-1$; let us denote this value by $\kappa$. By (\ref{eq-txd}), 
\beq \Pb A_{m,n} \rb = \sum _{\alpha = \pm 1} \Pb A_{m,n-1}; X_m(n) = \alpha 
\epsilon _n ;~\Delta S_{m,\kappa}^*= \alpha 2 \tilde X_{m-1}(\kappa + 1) \rb .
\label{eq-iidl} \eeq

The event $A_{m,n-1}$ can be written here as $B_{m,n-1} C_{m,n-1}$, 
where 
\[ B_{m,n-1} = \lb \tilde X_m(j) = \epsilon _j,~1 \le j \le T_m(\kappa ) 
 \rb , \]
\[ C_{m,n-1} = \lb X_m(j)= \alpha \epsilon _j,~T_m(\kappa ) + 1 \le j \le n-1 
 \rb . \] 
Definition (\ref{eq-txd}) shows that $B_{m,n-1}$ is determined by  
$\tilde X_{m-1}(j)$~ $(1 \le j \le \kappa )$ and $ X_m(j)$~ 
$(1 \le j \le T_m( \kappa ))$ the values of which
do not influence anything else in (\ref{eq-iidl}). 

Then we distinguish two cases according to the parity of $n$.

{\sc Case 1:} {\em $n$ is odd.} Then $n-1$ is even and $S_m(T_m(\kappa )) 
= S_m(n-1)$. Further, let  
$\tau _{m,r} = \min \lb j : j > 0, | S_m(r+j) - S_m(r) |  = 2 \rb $ and 
$\Delta S_m(r) = S_m(r + \tau _{m,r}) - S_m(r)$ for $r \ge 0$. 
Then $S_m(T_m(\kappa + 1)) = S_m(n-1+ \tau _{m,n-1})$ and 
$\Delta S_{m,\kappa }^* = \Delta S_m(n-1)$. These and the argument 
above shows that in (\ref{eq-iidl}) $A_{m,n-1}$ is independent of the other 
terms. Consequently, (\ref{eq-iidl}) simplifies as
\beq \Pb A_{m,n} \rb = 2 \Pb A_{m,n-1} \rb {1 \over 2} \sum_{\beta = \pm 1} 
 \Pb X_m(n) = \epsilon _n ;~\Delta S_m(n-1) = 2 \beta \rb , 
\label{eq-iids} \eeq
since the value of $\alpha $ is immaterial and $\Pb \tilde X_{m-1}( \kappa +1) 
= \beta \rb = 1/2$, independently of everything else here. 

Finally, (\ref{eq-ru2}) and (\ref{eq-ru3}) can be applied to (\ref{eq-iids}): 
\beqan  
\Pb A_{m,n} \rb   
& = & \Pb A_{m,n-1} \rb \sum_{\beta = \pm 1} \Pb \Delta 
 S_m(n-1) = 2 \beta \mid X_m(n) = \epsilon _n \rb ~\Pb X_m(n) = \epsilon _n 
 \rb  \\ 
& = & \Pb A_{m,n-1} \rb ({3 \over 4} + {1 \over 4}) {1 \over 2} 
 = {1 \over 2} \Pb A_{m,n-1} \rb , 
\eeqan  
independently of $\epsilon _n$. 

{\sc Case 2:} {\em $n$ is even.} Then $n-2$ is even and the argument in Case 
1 could be repeated with $n-2$ in place of $n-1$, with the only exception 
that in (\ref{eq-iidl}) we have an additional term $\tilde X_m(n-1) = 
\alpha X_m(n-1)$. Then instead of (\ref{eq-iids}) we arrive at 
\beqa \lefteqn 
{\Pb A_{m,n} \rb } \nonumber \\
& = & \Pb A_{m,n-2} \rb \sum_{\beta = \pm 1} \Pb X_m(n-1) = \epsilon _{n-1}, 
 X_m(n) = \epsilon _n ;~ \Delta S_m(n-2) = 2 \beta  \rb 
 \nonumber \\
& = & \Pb A_{m,n-2} \rb {1 \over 2^2} \sum_{\beta = \pm 1} \Pb \Delta 
 S_m(n-2) = 2 \beta \mid X_m(n-1) = \epsilon _{n-1}, X_m(n) = \epsilon _n 
 \rb . \label{eq-iid2} \eeqa

The conditional probability in (\ref{eq-iid2}) is 
\[ \begin{array}{llll}
1 & \mbox{if $~~\beta = \epsilon _{n-1} = \epsilon _n ;$} 
& 0 & \mbox{if $~~- \beta = \epsilon _{n-1} = \epsilon _n ;$} \\
1/2 & \mbox{if $~~\beta = \epsilon _{n-1} = - \epsilon _n ;$} 
& 1/2 & \mbox{if $~~- \beta = \epsilon _{n-1} = - \epsilon _n .$}  
\end{array} \]
Thus the sum in (\ref{eq-iid2}) becomes
\[ 1+0 = 1~~ \mbox{if $~~\epsilon _{n-1} = \epsilon _n;$}~~~~ 
 (1/2) + (1/2) = 1~~ \mbox{if $~~\epsilon _{n-1} = - \epsilon _n .$ } \]
In other words, the value of the sum in (\ref{eq-iid2}) 
is $1$, independently of $\epsilon _{n-1}$ and $\epsilon _n$. 

In summary, $\Pb A_{m,n} \rb = {1 \over 2} \Pb A_{m,n-1} \rb $ if $n$ is 
odd and $\Pb A_{m,n} \rb = {1 \over 4} \Pb A_{m,n-2} \rb $ if $n$ is even. 
Since $\Pb A_{m,0} \rb =1$, (\ref{eq-iid}) follows. $\Box$

%.................................

\vspace{10pt}

We mention that an other possibility to prove Lemma \ref{l-rw} is to 
introduce the random variables $Z_k={1 \over 2} \Delta S_{m,k-1}^* 
\tilde X_{m-1}(k)$ for $k \ge 1$. It can be shown that $Z_1, Z_2, \ldots $ is 
a sequence of independent and identically distributed random variables, 
$\Pb Z_k = 1 \rb = \Pb Z_k =-1 \rb = 1/2$, and this sequence is independent 
of the sequence $X_m(1), X_m(2), \ldots $ as well. Then we have 
$\tilde X_m(n) = Z_k X_m(n)$ for each $n$ such that $T_m(k-1) < n \le T_m(k)$ 
($k \ge 1$) and this implies (\ref{eq-iid}). 

The main property that was aimed when we introduced the ``twist'' method
easily follows from (\ref{eq-txd}) and (\ref{eq-tsmn}):
\beq \tilde S_m(T_m(k)) = \sum _{j=1}^k \tilde S_m(T_m(j)) - 
\tilde S_m(T_m(j-1)) = \sum _{j=1}^k 2 \tilde X_{m-1}(j)
= 2 \tilde S_{m-1}(k), 
\label{eq-refi} \eeq
for any $m \ge 1$ and $k \ge 0$.

\vspace{10 pt}

Now the second step of the approximation comes: ``shrinking''. As was 
discussed above, at the $m$th approximation the length of one step should be 
$1/2^m$ and the time needed for a step should be $1/2^{2m}$ (Figure~5). 
So we define the {\em $m$th approximation of the Wiener process} by 
\beq B_m\left( {t \over 2^{2m}}\right) = {1 \over 2^m} \tilde S_m(t)~~
 (t \ge 0, m \ge 0), \label{eq-bmt} \eeq
or equivalently, $B_m(t) = 2^{-m} \tilde S_m(t 2^{2m})$.
Basically, $B_m(t)$ is a model of Brownian motion on the set of points 
$x=j/2^m~(j \in \ZZ )$. 

Now (\ref{eq-refi}) becomes the following 
{\em refinement property\/}: 
\beq B_m\left( {T_m(k) \over 2^{2m}}\right) = B_{m-1}\left( {k \over 
2^{2(m-1)}}\right) , \label{eq-brefi} \eeq
for any $m \ge 1$ and $k \ge 0$.

\vspace{10 pt}

The remaining part of this section is devoted to showing the convergence
of the sequence $B_m(t)~(m=0,1,2,\ldots)$, and that the limiting process has 
the characterizing properties of the Wiener process. In proving these our 
basic tools will be some relatively simple, but powerful observations. 

First, often in the sequel the following crude, but still efficient 
estimate will be applied:
\beq \Pb \max_{1 \le j \le N} Z_j \ge t \rb = \Pb \bigcup _{j=1}^N \lb Z_j 
\ge t \rb \rb \le \sum _{j=1}^N \Pb Z_j \ge t \rb , \label{eq-crude} \eeq
which is valid for arbitrary random variables $Z_j$ and real number $t$.

The proofs of Lemmas \ref{l-tk} and \ref{l-bn} below essentially consist of 
the application of 
the following large deviation type estimate fulfilled by $S_n$ and 
$(T_k-4k)/ \sqrt 8$ according to Theorems \ref{thm-clt1}(b) and 
\ref{thm-clt2}(b). The previously mentioned exponential Chebyshev's 
inequalities (\ref{eq-chesn}) and (\ref{eq-chetk}) will be also used. 
Note that in the next lemma we have $a=2$ and $b=1$ for $S_n$ in 
(\ref{eq-chesn}) and $a=2$ and $b=1/2$ for $(T_k-4k)/ \sqrt 8$ in 
(\ref{eq-chetk}). 

%...............................

\begin{lemma} \label{l-ld}
Suppose that for $j \ge 0$, we have $\E (Z_j)=0$, $\Var (Z_j)=j$, and 
with some $a>0$ and $b>0$, 
\[ \Pb |Z_j| \ge t \rb \le 2 a^j e^{-bt} ~~ (t>0) \]
(exponential Chebyshev-type inequality). 

Assume as well that there exists a $j_0>0$ such that for any 
$j \ge j_0$, 
\[ \Pb |Z_j| / \sqrt j \ge x_j \rb \le e^{-x_j^2 /2}, \]
whenever $x_j \to \infty$ and $x_j=o(j^{1/6})$ as $j \to \infty$ (large  
deviation type inequality). 

Then for any $C>1$,
\beq \Pb \max _{0 \le j \le N} |Z_j| \ge  \sqrt{2CN \log N} \rb \le 
{2 \over N^{C-1}},. \label{eq-basle} \eeq
if $N$ is large enough, $N \ge N_0(C)$. 
\end{lemma}

{\sc Proof.}
The maximum in (\ref{eq-basle}) can be handled by the crude estimate 
(\ref{eq-crude}). Divide the resulting sum into two parts: one that can be 
estimated by a large deviation type inequality, and an other that will be 
estimated using exponential Chebyshev's inequality. For the large deviation 
part $x_j$ 
will be $\sqrt {2C \log N}$. Since $j \le N$, $j \to \infty $ implies that 
$N \to \infty $, and then $x_j \to \infty $ as well. If $j \ge \log ^4 N$, 
then the condition $x_j=o(j^{1/6})$ holds too, and the large deviation 
type inequality is applicable. Thus 
\beqan 
\lefteqn {\Pb \max _{0 \le j \le N} |Z_j| \ge \sqrt{2CN \log N} \rb} \\
 & \le & \sum _{j=0}^{\lfloor \log ^4 N \rfloor } 2 a^j \exp \left( -b \sqrt 
 {2CN \log N} \right) + \sum _{j= \lfloor \log ^4 N \rfloor }^N \Pb | Z_j|/ 
 \sqrt j \ge \sqrt {2C \log N} \rb \\
 & \le & {2a \over a-1} \exp \left( \log a~\log ^4 N - b \sqrt{2CN \log N} 
 \right) + N \exp (-C \log N) \le  2 ~ N^{1-C}
\eeqan
if $C>1$ and $N \ge N_0(C)$. ($\lfloor x \rfloor $ denotes the greatest 
integer less than or equal to $x$.) $\Box $

\vspace{10 pt}

%...............................

Note that the lemma and its proof are valid even when $N$ is not an integer. 
Here and afterwards we use the convention that if the upper limit of a 
sum is a real value $N$, then the sum goes until $\lfloor N \rfloor $. 

We mention that both inequalities among the assumptions of the previous Lemma
\ref{l-ld} hold for partial sums $Z_j$ of any sequence of independent and 
identically distributed random variables with expectation $0$, variance $1$ 
and a moment generating function which is finite 
for some $\pm u_0$. The fact that an exponential Chebyshev-type inequality 
should hold then can be seen from (\ref{eq-chee}) and (\ref{eq-mgfsn}), while 
the large deviation type estimate is shown to hold e.g. in 
\cite[Section XVI,6]{Feller 1966}.

The {\em first Borel--Cantelli lemma} \cite[Section VIII,3]{Feller 1968} 
is also an important tool, stating that if there is given an infinite 
sequence $A_1, A_2, \ldots $ of events such that 
$ \sum_{m=1}^{\infty } \Pb A_m \rb $ is finite, 
then {\em with probability $1$} only finitely many of the events occur. Or 
with an other widely used terminology: {\em almost surely} only finitely many 
of them will occur. 

\vspace{10 pt}

Now turning to the convergence proof, as the first step, it will be shown that 
the time instants $T_{m+1}(k)/2^{2(m+1)}$ will get arbitrarily 
close to the time instants $k/2^{2m} = 4k/2^{2(m+1)}$ as $m \to \infty $. 
By (\ref{eq-brefi}), 
this means that the next approximation not only visits the points $x=j/2^m~
(j \in \ZZ )$ in the same order, but the corresponding time instants will 
get arbitrarily close to each other as $m \to \infty $. Remember that by 
(\ref{eq-nbd}) and (\ref{eq-nbev}), $T_m(k)$ is the double of a negative 
binomial random variable, with expectation $4k$ and variance $8k$. Here Lemma 
\ref{l-ld} will be applied to $(T_m(k)-4k)/\sqrt 8$ with $N=K 2^{2m}$. So 
$\log N = \log K + (2 \log 2)m \le 1.5 m$ if $m$ is large enough, 
$m \ge m_0(K)$, and then $ \sqrt {2CN \log N} \le \sqrt {3CKm} \: 2^m$. 

%...............................

\begin{lemma} \label{l-tk}
(a) For any $C>1$, $K > 0$, and for any $m \ge m_0(C,K)$ we have
\beq \Pb \max _{0 \le k/2^{2m} \le K} |T_{m+1}(k) - 4k| \ge \sqrt {24 C K m} 
\: 2^m \rb < 2 \: (K 2^{2m})^{1-C}. \label{eq-pt} \eeq

(b) For any $K > 0$,
\beq \max _{0 \le k/2^{2m} \le K} \left| {T_{m+1}(k) \over 2^{2(m+1)}} - 
{k \over 2^{2m}} \right| <  \sqrt {2 K m} \: 2^{-m} 
\label{eq-ast} \eeq 
with probability $1$ for all but finitely many $m$. 
\end{lemma}

{\sc Proof.}

(a) (\ref{eq-pt}) is a direct consequence of Lemma \ref{l-ld}.

(b) Take for example $C=4/3$ in (a) and define the following events for 
$m \ge 0$:
\[ A_m = \lb \max _{0 \le k/2^{2m} \le K} |T_{m+1}(k) - 4k| \ge  \sqrt 
{32 K m} \: 2^m \rb . \]
By (\ref{eq-pt}), for $m \ge m_0(C,K)$, $\Pb A_m \rb < 2\: (K 2^{2m})^{-1/3}$.
Then $\sum _{m=0}^{\infty} \Pb A_m \rb < \infty$. 
Hence the Borel--Cantelli lemma implies that with probability $1$, only 
finitely many of the events $A_m$ occur. That is, almost surely for all but 
finitely many $m$ we have 
\[ \max _{0 \le k/2^{2m} \le K} |T_{m+1}(k) - 4k| <  \sqrt {32 K m} \: 2^m . 
\] 
This inequality is equivalent to (\ref{eq-ast}). $\Box$

%...............................

\vspace{10pt}

It seems to be important to emphasize a ``weakness'' of a statement like the 
one in Lemma \ref{l-tk}(b): we use the phrase ``all but finitely many $m$'' to 
indicate that the statement holds for 
every $m \ge m_0(\omega )$, where $m_0(\omega )$ may depend on the specific 
point $\omega $ of the sample space. 
In other words, one has no common, uniform lower bound for $m$ in general. 

Next we want to show that for any $j \ge 1$, $B_{n+j}(t)$ will be arbitrarily 
close to $B_n(t)$ as $n \to \infty $. Here again Lemma \ref{l-ld} will be 
applied, this time to a random walk $S_r$, with a properly 
chosen $N'$ and $C'$ (instead of $N=K 2^{2m}$ and $C$). Although the proof 
will be somewhat long, its basic idea is simple. Since  
$B_{m+1}(T_{m+1}(k)/2^{2(m+1)}) = B_m(k/2^{2m})$ by (\ref{eq-refi}), and 
the difference of the corresponding time instants here approaches zero fast  
as $m \to \infty $ by (\ref{eq-ast}), one can show that $B_m(t)$ and its 
refinement $B_{m+1}(t)$ will get very close to each other too. 

The following elementary fact that we need in the proof is discussed before 
stating the lemma: 
\beq \sum _{m=n}^{\infty } m 2^{-m/2} = (1/\sqrt 2) 
 \sum _{m=n}^{\infty } m \left( 1/\sqrt 2 \right) ^{m-1} < 4n 2^{-n/2}, 
 \label{eq-sum} \eeq
for $n \ge 15$. This can be shown by a routine application of power series:  
\[ \sum _{m=n}^{\infty } m x^{m-1} = {d \over dx} \sum _{m=n}^{\infty } x^m  
  = {d \over dx} \left( {x^n \over 1-x} \right) = n x^{n-1} \left( {1 \over 
  1-x} + {x \over n(1-x)^2} \right) . \]
Substituting $x=1/\sqrt 2$, one gets (\ref{eq-sum}) for $n \ge 15$. 

%...............................

\begin{lemma} \label{l-bn}
(a) For any $C \ge 3/2$, $K > 0$, and for any $n \ge n_0(C,K)$ we have
\beq \Pb \max _{0 \le k/2^{2n} \le K} |B_{n+1}(T_{n+1}(k)/2^{2(n+1)}) - 
 B_{n+1}(k/2^{2n})| \ge (1/8) n 2^{-n/2} \rb \le 3 (K 2^{2n})^{1-C} 
\label{eq-pbn} \eeq
and
\beq \Pb \max _{0 \le t \le K} |B_{n+j}(t) - B_n(t)| \ge  
n 2^{-n/2}~~\mbox{for some}~j \ge 1 \rb < 6 (K 2^{2n})^{1-C}. \label{eq-pb} 
\eeq

(b) For any $K > 0$, 
\beq \max _{0 \le t \le K} |B_{n+j}(t) - B_n(t)| 
 < n 2^{-n/2} , \label{eq-asb} \eeq 
with probability $1$ for all $j \ge 1$ and for all but finitely many $n$. 
\end{lemma}

{\sc Proof.} Let us consider first the difference between two consecutive 
approximations, $\max _{0 \le t \le K} |B_{m+1}(t) - B_m(t)|$. 
The maximum over real values $t$ can be approximated by the maximum over 
dyadic rational numbers $k/2^{2m}$. This is true because any sample-function 
$B_m(t; \omega )$ is a broken 
line such that, by (\ref{eq-bmt}), the magnitude of the increment between two 
consecutive points $k/2^{2m}$ and $(k+1)/2^{2m}$ is equal to $2^{-m}$.
Thus, taking the integer $t_m=\lfloor t 2^{2m} \rfloor $ for each 
$t \in [0,K]$, one has $t_m/2^{2m} \le t < (t_m+1)/2^{2m}$ and so 
$4t_m/2^{2(m+1)} \le t < (4t_m+4)/2^{2(m+1)}$. So we get
$|B_m(t)-B_m(t_m/2^{2m})| < 2^{-m}$ and 
$ \left| B_{m+1}(t)-B_{m+1}\left( 4t_m/2^{2(m+1)}\right) \right|$ 
$\le 4 \cdot 2^{-(m+1)}$ $= 2 \cdot 2^{-m}$. Hence 
\[  \max _{0 \le t \le K} |B_{m+1}(t) - B_m(t)| 
\le \max _{0 \le k/2^{2m} \le K} \left| B_{m+1}\left( 4k/2^{2(m+1)}\right) 
 - B_m\left( k/ 2^{2m}\right) \right| + 3 \cdot 2^{-m}. \]

Moreover, by (\ref{eq-brefi}) and (\ref{eq-bmt}) we have 
\beqa
 B_{m+1}\left(4k/2^{2(m+1)}\right) - B_m\left( k/2^{2m}\right) 
& = & B_{m+1}\left( 4k/2^{2(m+1)}\right) - B_{m+1}\left( T_{m+1}(k)/ 
 2^{2(m+1)}\right) \nonumber \\ 
& = & 2^{-(m+1)}\tilde S_{m+1}(4k) - 2^{-(m+1)}\tilde S_{m+1}(T_{m+1}(k)) . 
\label{eq-pb0} \eeqa

Thus 
\beqa \lefteqn
{\Pb \max _{0 \le t \le K} |B_{m+1}(t) - B_m(t)| \ge (1/4) m 2^{-m/2} \rb }
\nonumber \\
& \le & \Pb \max _{0 \le k/2^{2m} \le K} \left| B_{m+1}
 \left( 4k/2^{2(m+1)}\right) - B_m\left( k/2^{2m}\right) \right| 
 \ge (1/8) m 2^{-m/2} \rb \nonumber \\
& = & \Pb \max _{0 \le k \le K 2^{2m}} |\tilde S_{m+1}(4k) 
 - \tilde S_{m+1}(T_{m+1}(k))| \ge (1/4) m 2^{m/2} \rb 
\label{eq-pb1} \eeqa
if $m$ is large enough.  

By Lemma \ref{l-tk}, the probability of the event 
\[ A_m=\lb \max _{0 \le k \le K 2^{2m}}|T_{m+1}(k)-4k| \ge \sqrt{24CKm} 
\: 2^m \rb \] 
is very small for $m$ large. Therefore divide the last expression in 
(\ref{eq-pb1}) into two parts according to $A_m$ and $A_m^c$ 
(the complement of $A_m$): 
\beqa \lefteqn
{\Pb \max _{0 \le t \le K} |B_{m+1}(t) - B_m(t)| \ge (1/4) m 2^{-m/2} \rb }
 \nonumber \\
& \le & \Pb A_m^c ; \max _{0 \le k \le K 2^{2m}} \left| \tilde S_{m+1}(4k)
 - \tilde S_{m+1}(T_{m+1}(k)) \right| \ge (1/4) m 2^{m/2} \rb + \Pb A_m \rb
 \nonumber \\
& \le & \sum _{k=1}^{K 2^{2m}} \Pb \max _{\lb j:~|j-4k| 
 < \sqrt {24CKm} \: 2^m \rb } 
 \left| \tilde S_{m+1}(4k) - \tilde S_{m+1}(j) \right| \ge (1/4) m 2^{m/2} 
 \rb \nonumber \\
& & + 2 (K 2^{2m})^{1-C},
\label{eq-pb2} \eeqa
where the crude estimate (\ref{eq-crude}) and Lemma \ref{l-tk}(a) were used. 

Now apply Lemma \ref{l-ld} to $\tilde S_{m+1}(j)-\tilde S_{m+1}(4k)$ here,   
with suitably chosen $N'$ and $C'$. For $k$ fixed and $j>4k$, 
$\tilde S_{m+1}(j)-\tilde S_{m+1}(4k) = \sum _{r=4k+1}^j X_{m+1}(r)$ is a 
random walk of the form $S(j-4k)$. (The case $j<4k$ is symmetric.)
Since $|j-4k|<\sqrt{24CKm} \: 2^m$, $N'$ 
is taken as $\sqrt {24 CKm} \: 2^m$. (So $N'$ is roughly $\sqrt N$, where 
$N=K 2^{2m}$.) Then $\log N'=(1/2) \log (24CKm)$ $+ (\log 2) m \le m$ if  
$m$ is large enough, depending on $C$ and $K$. So 
\[ \sqrt {2C'N' \log N'} \le \sqrt {2C'm \sqrt {24CKm}\: 2^m} 
\le (1/4) m 2^{m/2}, \] 
if $m$ is large enough, depending on $C$, $C'$, and $K$. Then it follows by 
Lemma \ref{l-ld} that 
\beq \Pb \max _{0 \le r \le \sqrt {24 CKm} \: 2^m} |S(r)| \ge (1/4) m 2^{m/2}
 \rb \le 2 (\sqrt {24CKm} \: 2^m)^{1-C'}. \label{eq-sqrtn} \eeq

The second term of the error probability in (\ref{eq-pb2}) is $2 (K 2^{2m})
^{1-C}=2 N^{1-C}$, while (\ref{eq-sqrtn}) indicates that the first term is at  
most $K 2^{2m}\cdot 2 \cdot 2 (\sqrt {24CKm} \: 2^m)^{1-C'}$ $\le N (\sqrt N) 
^{1-C'}$ if 
$C'>1$ and $m$ is large enough. To make the two error terms to be of the same 
order, choose $1+(1-C')/2 = 1-C$, i.e. $C'=2C+1$. Thus 
(\ref{eq-pb2}) becomes
\[ \Pb \max _{0 \le t \le K} |B_{m+1}(t) - B_m(t)| \ge (1/4) m 2^{-m/2} \rb
\le 3 (K 2^{2m})^{1-C} , \]
for any $m$ large enough, depending on $C$ and $K$. Comparing this to 
(\ref{eq-pb0}) and (\ref{eq-pb1}) one obtains (\ref{eq-pbn}).

By (\ref{eq-sum}), $ \max _{0 \le t \le K} |B_{m+1}(t) - B_m(t)| 
< (1/4) m 2^{-m/2}$ for all $m \ge n \ge 15$ would imply that 
\beqan
& & \max _{0 \le t \le K} |B_{n+j}(t) - B_n(t)| = \max _{0 \le t \le K} 
 \left| \sum _{m=n}^{n+j-1} B_{m+1}(t) - B_m(t) \right| \\
& \le & \sum _{m=n}^{n+j-1} \max _{0 \le t \le K} |B_{m+1}(t) - B_m(t)| 
 < \sum _{m=n}^{\infty } (1/4) m 2^{-m/2} 
 < n 2^{-n/2}, 
\eeqan 
for any $j \ge 1$. So we conclude that  
\beqan \lefteqn
{\Pb \max _{0 \le t \le K} |B_{n+j}(t) - B_n(t)| \ge n 2^{-n/2}
 ~~\mbox{for some}~j \ge 1 \rb } \\
& \le & \sum _{m=n}^{ \infty }  
 \Pb \max _{0 \le t \le K} |B_{m+1}(t) - B_m(t)| \ge (1/4) m 2^{-m/2} \rb \\
& \le & \sum _{m=n}^{\infty } 3 (K 2^{2m})^{1-C} = 3 (K 2^{2n})^{1-C} 
{1 \over 1-2^{2(1-C)}} < 6 (K 2^{2n})^{1-C} 
\eeqan
if $C \ge 3/2$ (say), for any $n \ge n_0(C,K)$. This proves (\ref{eq-pb}). 

The statement in (b) follows from (\ref{eq-pb}) by the Borel--Cantelli 
lemma, as in the proof of Lemma \ref{l-tk}. $\Box$

%............................

\vspace{10pt}

Now we are ready to establish the existence of the Wiener process, which is a 
continuous model of Brownian motion. An important consequence of 
(\ref{eq-asb}) is that the difference between the Wiener process and $B_n(t)$ 
is smaller than a constant multiple of $\log N /\sqrt[4] N$, where 
$N=K 2^{2n}$, see (\ref{eq-asw}) below. 

%............................

\begin{thm} \label{thm-w}
As $n \to \infty $, with probability $1$ (that is, for almost all $\omega \in 
\Omega $) and for all $t \in [0, \infty)$ the sample-functions 
$B_n(t;\omega )$ converge to a sample-function $W(t;\omega)$ such that 

(i) $W(0;\omega )=0$, $W(t;\omega )$ is a continuous function of $t$ on the 
interval $[0, \infty )$;

(ii) for any $0 \le s < t$, $W(t)-W(s)$ is a normally distributed random 
variable with expectation $0$ and variance $t-s$;

(iii) for any $0 \le s < t \le u < v$, the increments $W(t)-W(s)$ and 
$W(v)-W(u)$ are independent random variables. 

By definition, $W(t)~(t \ge 0)$ is called the {\bf Wiener process}.

Further, we have the following estimates for the difference of the Wiener 
process and its approximations.

(a) For any $C \ge 3/2$, $K > 0$, and for any $n \ge n_0(C,K)$ we have
\beq \Pb \max _{0 \le t \le K} |W(t) - B_n(t)| \ge n 2^{-n/2} \rb \le 
6 (K 2^{2n})^{1-C} . \label{eq-pw} \eeq

(b) For any $K > 0$,
\beq \max _{0 \le t \le K} |W(t) - B_n(t)| \le n 2^{-n/2} , 
\label{eq-asw} \eeq 
with probability $1$ for all but finitely many $n$. 
\end{thm}

{\sc Proof.}
Lemma \ref{l-bn}(b) shows that for almost all $\omega \in \Omega $, the 
sequence $B_n(t; \omega )$ converges for any $t \ge 0$ as $n \to \infty $. Let 
us denote the limit by $W(t; \omega )$. On a probability zero $\omega $-set 
the limit possibly does not exist, there one can define $W(t; \omega ) = 0$ 
for any $t \ge 0$. Since $B_n(0; \omega )=0$ for any $n$, it follows that 
$W(0; \omega )=0$ for any $\omega \in \Omega $.

Taking $j \to \infty $ in (\ref{eq-pb}), (\ref{eq-pw}) follows. By
(\ref{eq-asb}), the convergence of $B_n(t)$ is uniform on any bounded interval 
$[0, K]$, more exactly, for any $K > 0$ we have (\ref{eq-asw}) with 
probability $1$. Textbooks on 
advanced calculus, like W. Rudin's \cite[Section 7.12]{Rudin 1977} show that 
the limit function of a uniformly convergent sequence of continuous functions 
is also continuous. This proves (i). 

Now we turn to the proof of (ii). Take arbitrary $t > s \ge 0$ and $x$ real. 
With $K>t$ fixed, (\ref{eq-pw}) shows that for any  
$\delta > 0$ there exists an $n \ge n_0(C,K)$ such that 
\beq \Pb \max _{0 \le u \le K} |W(u)-B_n(u)| \ge \delta \rb < \delta . 
\label{eq-wb1} \eeq
Since 
\[ \Pb W(t)-W(s) \le x \rb = \Pb B_n(t)-B_n(s) \le x-(W(t)-B_n(t))
  +(W(s)-B_n(s)) \rb , \]
(\ref{eq-wb1}) implies that 
\beqa 
      \Pb B_n(t)-B_n(s) \le x-2 \delta \rb -2 \delta 
& \le & \Pb W(t)-W(s) \le x \rb \nonumber \\
& \le & \Pb B_n(t)-B_n(s) \le x+2 \delta \rb +2 \delta . 
\label{eq-wb2} \eeqa
This indicates that the distribution function of $W(t)-W(s)$ can be 
eventually obtained from the distribution function of 
\beq B_n(t)-B_n(s)= 2^{-n} \tilde S_n(2^{2n}t)-2^{-n} \tilde S_n(2^{2n}s).
\label{eq-wb3} \eeq 

Take the non-negative integers $j_n= \lfloor 2^{2n}t \rfloor $ and 
$i_n= \lfloor 2^{2n}s \rfloor $, $j_n \ge i_n$. Then (\ref{eq-wb3}) differs 
from 
\beq 2^{-n}(\tilde S_n(j_n)- \tilde S_n(i_n)) = 2^{-n}
 \sum _{k=i_n+1} ^{j_n} \tilde X_k \label{eq-wb4} \eeq
by an error not more than $2\cdot 2^{-n} < \delta $. (We can assume that $n$ 
was chosen so.) Also, $j_n-i_n$ differs from $2^{2n}(t-s)$ by at most $1$. 
In particular, $j_n-i_n \to \infty$ as $n \to \infty$.  

If $n$ is large enough (we can assume again that $n$ was chosen so), 
by Theorem \ref{thm-clt1}(a), for any fixed real $x'$ we have 
\beq \Phi (x') - \delta \le \Pb {1 \over \sqrt {j_n-i_n}} \sum _{k=i_n+1} 
^{j_n} \tilde X_k \le x' \rb \le \Phi (x') + \delta . 
\label{eq-wb5} \eeq
Here $\sqrt {j_n-i_n}$ can be approximated by $2^n\sqrt {t-s}$ if $n$ is 
large enough, 
\beq 1-\delta < \sqrt {j_n-i_n-1 \over j_n-i_n} \le {2^n \sqrt {t-s} \over
 \sqrt {j_n-i_n}} \le \sqrt {j_n-i_n+1 \over j_n-i_n} <1+\delta . 
\label{eq-wb6} \eeq 

Combining formulae (\ref{eq-wb3})-(\ref{eq-wb6}) we obtain that 
\[ \Phi \left( (1-\delta ){x \over \sqrt{t-s}}-\delta \right) - \delta \le 
\Pb B_n(t)-B_n(s) \le x \rb \le \Phi \left( (1+\delta ){x \over \sqrt{t-s}}
+\delta \right) + \delta . 
\]
This shows that the distribution of $B_n(t)-B_n(s)$ is asymptotically normal 
with mean $0$ and variance $t-s$ as $n \to \infty $. Moreover, by 
(\ref{eq-wb2}), the distribution of $W(t)-W(s)$ is {\em exactly} normal 
with mean $0$ and variance $t-s$, since $\delta $ can be made arbitrarily 
small if $n$ is large enough:
\[ \Pb W(t)-W(s) \le x \rb = \Phi \left( {x \over \sqrt{t-s}} \right) . \]
This proves (ii). 

Finally, (iii) can be proved similarly as (ii) above. Taking arbitrary 
$v>u \ge t>s \ge 0$ and $x$, $y$ real numbers, 
\beq \Pb W(t)-W(s) \le x,~W(v)-W(u) \le y \rb \label{eq-wb7} \eeq
can be approximated by a probability of the form 
\[ \Pb B_n(t)-B_n(s) \le x~,~B_n(v)-B_n(u) \le y \rb \]
arbitrarily well if $n$ is large enough, just like in (\ref{eq-wb2}). In 
turn, like in (\ref{eq-wb4}), the latter can be estimated arbitrarily well 
by a probability of the form 
\beq \Pb {1 \over \sqrt {j_n-i_n}} \sum _{k=i_n+1}^{j_n} \tilde X_k \le x'~,
 ~{1 \over \sqrt {r_n-q_n}} \sum _{k=q_n+1}^{r_n} \tilde X_k \le y' \rb , 
\label{eq-wb8} \eeq
where $i_n= \lfloor 2^{2n}s \rfloor \le j_n= \lfloor 2^{2n}t \rfloor \le 
q_n= \lfloor 2^{2n}u \rfloor \le r_n= \lfloor 2^{2n}v \rfloor $. 

Since there are no common terms in the first and the second sum of 
(\ref{eq-wb8}), the two sums are independent. Thus (\ref{eq-wb8}) is 
equal to 
\[ \Pb {1 \over \sqrt {j_n-i_n}} \sum _{k=i_n+1}^{j_n} \tilde X_k \le x' \rb 
~\cdot ~\Pb {1 \over \sqrt {r_n-q_n}} \sum _{k=q_n+1}^{r_n} \tilde X_k \le y' 
\rb ,\]
which can be made arbitrarily close to 
\beq \Pb W(t)-W(s) \le x \rb~\cdot ~\Pb W(v)-W(u) \le y \rb. 
\label{eq-wb9} \eeq
Since errors in the approximations can be made arbitrarily small, 
(\ref{eq-wb7}) and (\ref{eq-wb9}) must agree for any real $x$ and $y$. This 
proves (iii). $\Box $

%...............................

\vspace{10pt}

Note that properties (ii) and (iii) are often rephrased in the way that the 
Wiener-process is a Gaussian process with independent and stationary 
increments. It can be proved \cite[Section 1.5]{Knight 1981} that properties 
(i), (ii), and (iii) characterize the Wiener process. In other words, any 
construction to the Wiener process gives essentially the same process that 
was constructed above. 

%-------------------------------------------------------------------------

\section{From the Wiener Process to Random Walks}

Now we are going to check whether the Wiener process as a model of Brownian 
motion has the properties described in the introduction of Section 
\ref{se-rwtw}. Namely, we would want to find the sequence of shrunk random 
walks $B_m(k 2^{-2m})$ in $W(t)$. 

Let $s(1)$ be the first (random) time instant where the magnitude of the 
Wiener process is $1$: $s(1) = \min \lb s>0: |W(s)| = 1 \rb .$  
The continuity and increment characteristics of the Wiener process imply that 
$s(1)$ exists with probability $1$. Clearly, each shrunk random walk $B_m(t)$ 
has the symmetry property that reflecting all its sample-functions to the 
time axis, one gets the same process. $W(t)$ as a limiting process of shrunk  
random walks inherits this feature. Therefore setting $X(1)=W(s(1))$, 
$\Pb X(1)=1 \rb = \Pb X(1)=-1 \rb = 1/2.$

Inductively, starting with $s(0)=0$, if $s(k-1)$ is given, define 
the random time instant
\[ s(k) = \min \lb s: s>s(k-1), |W(s)-W(s(k-1))|=1 \rb ~~~(k \ge 1). \]
As above, $s(k)$ exists with probability $1$. Setting $X(k)=W(s(k))
-W(s(k-1))$, it is heuristically clear that $\Pb X(k)=1 \rb = \Pb X(k)=-1 \rb 
= 1/2$, and $X(k)$ is independent of $X(1),X(2), \ldots ,X(k-1)$.  

This way one gets a random walk $S(k)=W(s(k))=X(1)+X(2)+ \cdots X(k)$ 
$(k \ge 0)$ from the Wiener process. Using a more technical phrase, 
by this method, based on {\em first passage times}, one can {\em imbed} a 
random walk into the Wiener process;
it is a special case of the famous Skorohod imbedding, see e.g. 
\cite[Section 13.3]{Breiman 1968}. 

Quite similarly, one can imbed $B_m(k2^{-2m})$ into $W(t)$ for any $m \ge 0$ 
by setting $s_m(0)=0$, 
\beq s_m(k) = \min \lb s: s>s_m(k-1), |W(s)-W(s_m(k-1))|=2^{-m} \rb (k \ge 1), 
\label{eq-smk} \eeq 
and $B_m(k2^{-2m})=W(s_m(k))$ $(k \ge 0)$. 

However, instead of proving all necessary details about Skorohod 
imbedding briefly described above, we will define an other imbedding method 
better suited to our approach. It will turn out that our imbedding is  
essentially equivalent to the Skorohod imbedding. 

Our task requires a more careful analysis of the waiting times $T_m(k)$ first. 
Recall the refinement property (\ref{eq-brefi}) of $B_m(t)$. Continuing that, 
we get  
\beqa 
B_m(k 2^{-2m}) & = & B_{m+1}(2^{-2(m+1)} T_{m+1}(k)) 
= B_{m+2}(2^{-2(m+2)} T_{m+2}(T_{m+1}(k))) = \cdots  
 \nonumber\\
& = & B_n(2^{-2n} T_n(T_{n-1}(\cdots (T_{m+1}(k)) \cdots ))) ,
\label{eq-bsame} \eeqa
where $k\ge 0$ and $n>m \ge 0$. In other words, $B_n(t)$, $n>m$, visits the 
same dyadic points $k 2^{-m}$ in the same order as $B_m(t)$, only the 
corresponding time instants can differ. 

To simplify the notation, let us introduce  
\[ T_{m,n}(k) = T_n(T_{n-1}(\cdots (T_{m+1}(k)) \cdots ))~~~(n>m \ge 0,~
 k \ge 0) . \]
Then (\ref{eq-bsame}) becomes 
\beq B_m(k 2^{-2m}) = B_n(2^{-2n} T_{m,n}(k)) ~~~(n>m \ge 0, ~ k \ge 0). 
\label{eq-bsamec} \eeq
The next lemma considers {\em time lags} of the form $2^{-2n} T_{m,n}(k)
-k 2^{-2m}$. 

Note that in the proofs of the next two lemmas we make use of the following 
simple inequality, valid for arbitrary random variables $Z_j$, real numbers 
$t_j$, and events $A_j=\lb Z_j \ge t_j \rb $:
\beqa 
& & \Pb Z_j \ge t_j ~~\mbox{for some $j \ge 1$} \rb   
 = \Pb \bigcup _{j=1}^{\infty } A_j \rb \nonumber \\ 
& = & \Pb A_1 \rb + \Pb A_1^c A_2 \rb + \cdots 
 + \Pb A_1^c \cdots  A_j^c A_{j+1} \rb + \cdots \nonumber \\
& \le & \Pb Z_1 \ge t_1 \rb 
+ \sum _{j=1}^{\infty } \Pb  Z_j < t_j, Z_{j+1} \ge t_{j+1} \rb . 
\label{eq-pineq} \eeqa

%...........................

\begin{lemma} \label{l-tmn}
(a) For any $C \ge 3/2$, $K>0$, and $m \ge 0$ take the following subset of 
the sample space:   
\beq A_m = \lb \max _{0 \le k 2^{-2m} \le K} |2^{-2n} T_{m,n}(k) 
 - k 2^{-2m}| < \sqrt {18CKm} \: 2^{-m} ~~\mbox{for all $n>m$} \rb . 
\label{eq-am} \eeq 
Then for any $m \ge m_0(C,K)$ we have 
\beq \Pb A_m \rb \ge 1 - 4 (K 2^{2m})^{1-C}. \label{eq-tmn} \eeq

(b) For any $K>0$, 
\[ \max _{0 \le k 2^{-2m} \le K} |2^{-2n} T_{m,n}(k) - k 2^{-2m}| <
\sqrt {27Km} \: 2^{-m} \]
with probability $1$ for all $n>m$, for all but finitely many $m$. 
\end{lemma}

{\sc Proof.} 
Take any $\beta $, $1/2 < \beta < 1$, and $K'>K$; for example $\beta = 2/3$ 
and $K'=(4/3)K$ will do. Set $\alpha _j = 1 + \beta + \beta ^2 + \cdots + 
\beta ^j$ for $j \ge 0$,
\[ Z_n = \max _{0 \le k2^{-2m} \le K} |T_{m,n}(k) - k 2^{2(n-m)}|,~~~
    t_n = \alpha _{n-m-1} \sqrt {24CK'm} \: 2^{2n-m-2}, \]
and $Y_{n+1} = \max _{0 \le k2^{-2m} \le K} |T_{n+1}(T_{m,n}(k)) - 4 
T_{m,n}(k)|$ for $n>m \ge 0$. 

First, by Lemma \ref{l-tk}(a), 
\[ \Pb Z_{m+1} \ge t_{m+1} \rb = \Pb \max _{0 \le k2^{-2m} \le K} |T_{m+1}(k) - 
 4k| \ge \sqrt {24CK'm} \: 2^m \rb \le 2(K2^{2m})^{1-C} \]
if $m$ is large enough. 

Second, by the triangle inequality, $Z_{n+1} \le 4 Z_n + Y_{n+1}$ for any 
$n>m$. So 
\[ \Pb Z_n < t_n, Z_{n+1} \ge t_{n+1} \rb 
\le \Pb Z_n < t_n, Y_{n+1} \ge t_{n+1} - 4t_n \rb . \]
If $Z_n < t_n$, then setting $j=T_{m,n}(k)$, 
\beqan
j2^{-2n} & < & 2^{-2n}(k2^{2(n-m)} + t_n) 
 = k2^{-2m} + \alpha _{n-m-1} \sqrt {24CK'm} \: 2^{-m-2} \\
& \le & K + 3 \sqrt{2CKm} \: 2^{-m} \le (4/3)K = K'
\eeqan
holds for $m \ge m_0(C,K)$, since $\alpha _r < 1/(1- \beta ) = 3$ (if $\beta 
= 2/3$) for any $r \ge 0$. 
Applying these first and Lemma \ref{l-tk}(a) last, for $n>m \ge m_0(C,K)$ we 
get that  
\beqan \lefteqn
{\Pb Z_n < t_n, Z_{n+1} \ge t_{n+1} \rb} \\ 
& \le & \Pb \max _{0 \le j2^{-2n} \le K'} |T_{n+1}(j) - 4 j| \ge 
 \sqrt {24CK'm} \: 2^n 2^{n-m} (\alpha _{n-m} - \alpha _{n-m-1}) \rb \\
& \le & \Pb \max _{0 \le j \le K' 2^{2n}} |T_{n+1}(j) - 4 j| \ge 
\sqrt {24CK'n} \: 2^n \rb \le 2 (K'2^{2n})^{1-C}. 
\eeqan
In the second inequality above we used that $\sqrt m 2^{n-m} (\alpha _{n-m} - 
\alpha _{n-m-1}) = (2 \beta )^{n-m} \sqrt m \ge \sqrt n $, 
which follows from the inequality  
$(2 \beta )^{n-m} = (4/3)^{n-m} \ge \sqrt {1+(n-m)/m}$, valid for any 
$n>m \ge 2$ (if $\beta = 2/3$). 

Combining the results above, 
\beqan \lefteqn
{\Pb \max _{0 \le k 2^{-2m} \le K} |2^{-2n} T_{m,n}(k) - k 2^{-2m}| \ge 
\sqrt {18CKm} \: 2^{-m} ~~\mbox{for some $n>m$} \rb } \\
& = & \Pb \max _{0 \le k 2^{-2m} \le K} |T_{m,n}(k) - k 2^{2(n-m)}| \ge 
 3 \sqrt {24CK'm} \: 2^{2n-m-2} ~~\mbox{for some $n>m$} \rb \\
& \le & \Pb Z_n \ge t_n ~~\mbox{for some $n \ge m+1$} \rb \\
& \le & \Pb Z_{m+1} \ge t_{m+1} \rb 
 + \sum _{n=m+1}^{\infty } \Pb  Z_n < t_n, Z_{n+1} \ge t_{n+1} \rb \\
& \le & \sum _{n=m}^{\infty } 2 (K2^{2n})^{1-C} = 2 (K2^{2m})^{1-C} 
 {1 \over 1 - 2^{2(1-C)}} \le 4 (K 2^{2m})^{1-C} 
\eeqan
if $C \ge 3/2 $, say. This proves (a). 

The statement in (b) follows by the Borel--Cantelli lemma. $\Box $

%.............................

\vspace{10pt}

As (\ref{eq-bsamec}) shows, $B_n(2^{-2n} T_{m,n}(k))=B_m(k 2^{-2m})$ for any
$k \ge 0$ and for any $n>m\ge 0$. A natural question, important particularly 
when looking at increments during short time intervals, is that how much time 
it takes for $B_n(t)$ to go from the point $B_m((k-1) 2^{-2m})$ to its 
next ``$m$-neighbor' $B_m(k 2^{-2m})$. Is this time significantly different 
from $2^{-2m}$ for large values of $m$? Introducing the notation 
\beq \tau _{m,n}(k)=T_{m,n}(k) - T_{m,n}(k-1) ~~~(k \ge 1,~n>m\ge 0), 
\label{eq-delk} \eeq
{\em the $n$th time differences of the $m$-neighbors} are $2^{-2n} 
\tau _{m,n}(k)$ $(k \ge 1)$. 
Note that $T_{m,n}(k) = \sum _{j=1}^k \tau _{m,n}(j)$, where, as can be 
seen from the construction and the argument below, the terms are independent 
and have the same distribution.

Let us look at $ \tau _{m,n}(k)$ more closely. If $n=m+1$, 
\beq \tau _{m,m+1}(k) = T_{m+1}(k) - T_{m+1}(k-1) = \tau _{m+1}(k), 
\label{eq-taum1} \eeq  
which is the double of a geometric random variable 
with parameter $p=1/2$, see (\ref{eq-geom}). That is, $2^{-2(m+1)} \tau _{m+1}
(k)$ is the length of the time period that corresponds to the time interval 
$[(k-1) 2^{-2m},~ k 2^{-2m}]$ after the next refinement of the construction. 

Similarly, each unit in $\tau _{m+1}(k)$ will bring some $\tau _{m+2}(r)$ 
``offsprings'' after the following refinement, and so on. Hence if 
$n>m$ is arbitrary, then given $T_{m,n}(k-1) = j$ for some integer $j \ge 0$, 
we have  
\beq \tau _{m,n+1}(k) = T_{n+1}(j+ \tau _{m,n}(k)) - T_{n+1}(j) 
 = \sum _{r=1}^{ \tau _{m,n}(k)} \tau _{n+1}(j+r). \label{eq-taumn} \eeq
For given $ \tau _{m,n}(k)=s$ ($s > 0$, even) its conditional distribution 
is the same as the distribution of a random variable $T_{s}$ which is the 
double of a negative binomial random variable with parameters 
$s$ and $p=1/2$, described by (\ref{eq-nebi}) and (\ref{eq-nbd}). Note 
that this conditional distribution of $\tau _{m,n+1}(k)$ is independent of 
the value of $T_{m,n}(k-1)$. 

Though we will not explicitly use them, it is instructive to determine 
some further properties of a ``prototype'' $\tau _{m,n}=\tau _{m,n}(1)$. 
A recursive formula can be given for its expectation by (\ref{eq-taumn}) 
and the full expectation formula: 
\[ \mu _{n+1} = \E (\tau _{m,n+1}) = \E (\E (\tau _{m,n+1}~|~\tau _{m,n})) 
 = \E (\tau _{m,n} \E (\tau _{n+1}(r))) = 4 \mu _n . \] 
Since $\mu _{m+1} = \E (\tau _{m+1}(r)) = 4$, it follows that 
\[ \mu _n = \E (\tau _{m,n}) = 2^{2(n-m)}. \] 

This argument also implies that 
\[ \E (2^{-2(n+1)} \tau _{m,n+1} ~|~ 2^{-2n} \tau _{m,n}) = 2^{-2n} 
 \tau _{m,n}. \] 
These show that the sequence $\left( 2^{-2n} \tau _{m,n} \right) _{n=m+1}
^{\infty }$ is a so-called {\em martingale}. Therefore a famous 
martingale convergence theorem \cite[Section 5.4]{Breiman 1968} implies 
that this sequence converges to a random variable $t_m$ as $n \to \infty $, 
with probability $1$, and $t_m$ has finite expectation. 

We mention that a similar recursion can be obtained for the variance that 
results 
\[ \Var ( 2^{-2n} \tau _{m,n}) < {2 \over 3} 2^{-4m} . \]  

The next lemma gives an upper bound for the $n$th time differences of the 
$m$-neighbors by showing that during arbitrary many refinements, 
they cannot be ``much'' larger than $h=2^{-2m}$, the original time difference 
of the $m$-neighbors. More accurately, they are less than a multiple of 
$h^{1-\delta }$, where $\delta >0$ arbitrary.

%............................

\begin{lemma} \label{l-delta}
(a) For any $K>0$, $\delta $ such that $0 < \delta <1$, and $C>2/ \delta $ 
we have  
\[ \Pb \max _{1 \le k 2^{-2m} \le K} |2^{-2n}\tau _{m,n}(k) -  2^{-2m}| 
 \ge 3C 2^{-2m(1- \delta )} ~~~\mbox{for some $n>m$} \rb 
 \le {K \over 10} 2^{-2m( \delta C - 2)} . \]
(b) For any $K>0$, and $\delta $ such that $0 < \delta <1$, 
\[ \max _{1 \le k 2^{-2m} \le K} |2^{-2n} \tau _{m,n}(k) -  2^{-2m}| 
 < {7 \over \delta } 2^{-2m(1- \delta )}, \] 
with probability $1$ for all $n>m$, for all but finitely many $m$.
\end{lemma}

{\sc Proof.}
This proof is very similar to the proof of Lemma \ref{l-tmn}. 
Take any $\beta $, $1/2 < \beta < 1$; for example $\beta = 2/3$ will do. 
Set $\alpha _j = 1 + \beta + \beta ^2 + \cdots + \beta ^j$ for $j \ge 0$. 
For any $m \ge 0$, consider an arbitrary $k$, $1 \le k \le K 2^{2m}$. (The 
distribution of $\tau _{m,n}(k)$ does not depend on $k$.) Let 
\[ Z_n = |\tau _{m,n}(k) -  2^{2(n-m)}|,~~~
    t_n = \alpha _{n-m-1} C 2^{2 \delta m} 2^{2(n-m)}, \]
and $Y_{n+1} = |\tau _{m,n+1}(k)) - 4 \tau _{m,n}(k)|$ for $n>m \ge 0$. 
We want to apply inequality (\ref{eq-pineq}). 

First take $n=m+1$. By (\ref{eq-taum1}), ${1 \over 2} \tau _{m,m+1}(k)$ is 
a geometric random variable with parameter $p=1/2$. Then  
\beqan 
\Pb Z_{m+1} \ge t_{m+1} \rb & = & \Pb |\tau _{m,m+1}(k) -4| \ge 4C 
 2^{2 \delta m} \rb \\
& = & \Pb {1 \over 2} \tau _{m+1}(k) \ge 2 + 2C 2^{2 \delta m} \rb < 
 2^{-4C \delta m} ,  
\eeqan
because of the basic property of the tail of a geometric distribution. 

Second, let $n>m$ be arbitrary. By the triangle inequality, $Z_{n+1} \le 
4 Z_n + Y_{n+1}$. So we obtain 
\beqa \lefteqn 
{\Pb Z_n < t_n, Z_{n+1} \ge t_{n+1} \rb \le \Pb Z_n < t_n, Y_{n+1} \ge 
t_{n+1}-4t_n \rb } \nonumber \\
& \le & \sum _{s=1}^{2^{2(n-m)}+t_n} \Pb |T_s - 4s| \ge t_{n+1}-4t_n ~|~ 
\tau _{m,n}(k)=s \rb ~ \Pb \tau _{m,n}(k)=s \rb \nonumber \\
& \le & \Pb \max _{1 \le s \le 2^{2(n-m)}+t_n} |T_s-4s|/ \sqrt 8 \ge 
(t_{n+1}-4t_n)/ \sqrt 8 \rb ~ \sum _{s=1}^{\infty } \Pb \tau _{m,n}(k)=s \rb ,
\label{eq-taumn1} \eeqa
where we applied (\ref{eq-taumn}) and the conditional distribution of 
$\tau _{m,n+1}(k)$ mentioned there.  

The sum in (\ref{eq-taumn1}) is $1$. Therefore we want to estimate the 
probability of the maximum there by using Lemma \ref{l-ld} with 
$N=4C 2^{2 \delta m} 2^{2(n-m)}$. This $N$ is larger than $2^{2(n-m)}+t_n$ 
if $m$ is large enough, depending on $\delta $. (Remember that 
$\alpha _j < 3$ for any $j \ge 0$ if $\beta =2/3$.) To apply Lemma 
\ref{l-ld} we have to compare $\sqrt {2CN \log N}$ to the right hand side 
of the inequality in (\ref{eq-taumn1}): 
\[ {\sqrt {2CN \log N} \over (t_{n+1}-4t_n)/\sqrt 8} = 2 \left( {(n-m) 
\log 4 + m \delta \log 4 + \log(4C) \over (4/3)^{2(n-m)} 2^{2 \delta m}} 
\right)^{1/2}, \] 
which is less than $1$ for all $n>m$ if $m$ is large enough, depending on 
$\delta $. 

Thus Lemma \ref{l-ld} gives that 
\beqan \Pb Z_n < t_n, Z_{n+1} \ge t_{n+1} \rb & \le & \Pb \max _{1 \le s 
\le N} |T_s-4s|/ \sqrt 8 \ge \sqrt {2CN \log N} \rb \\ 
& \le & 2 N^{1-C} =  2 (4C 2^{2 \delta m} 2^{2(n-m)})^{1-C} \eeqan
for all $n>m$ as $m \ge m_0(\delta , C)$. 

Summing up for $n>m$, we obtain the following estimate for any given $k$, 
$1 \le k 2^{-2m} \le K$, as $m \ge m_0(\delta , C)$: 
\beqan \lefteqn
{\Pb  |2^{-2n} \tau _{m,n}(k) - 2^{-2m}| \ge 3C 2^{-2m(1-\delta )} 
 ~~~\mbox{for some $n>m$} \rb } \\
& \le & \Pb Z_n \ge t_n ~~\mbox{for some $n \ge m+1$} \rb \\
& \le & \Pb Z_{m+1} \ge t_{m+1} \rb 
 + \sum _{n=m+1}^{\infty } \Pb  Z_n < t_n, Z_{n+1} \ge t_{n+1} \rb \\
& \le & 2^{-4C \delta m} + 2^{2 \delta (1-C) m} 2 (4C)^{1-C} \sum _{n=m+1}
^{\infty } 2^{2(1-C)(n-m)} \\  
& < & (1/10)  2^{-2m(\delta C - 1)} , 
\eeqan
where we took into consideration that $0<\delta <1$, $C>2$, $\delta C>2$ 
by our assumptions. 

Finally, the statement in (a) can be obtained by an application of the 
crude inequality (\ref{eq-crude}), 
\[ \Pb \max _{1 \le k \le K 2^{2m}} |2^{-2n}\tau _{m,n}(k) -  2^{-2m}| 
 \ge 3C 2^{-2m(1- \delta )} ~~~\mbox{for some $n>m$} \rb 
 \le K 2^{2m} {1 \over 10} 2^{-2m(\delta C - 1)} , \]
as $m \ge m_0(\delta, C)$, which is equivalent to (a). 

The statement in (b) follows by the Borel--Cantelli lemma with $C=7/(3 
\delta )$ (say). $\Box $

%.............................

\vspace{10pt}

Now we define a certain {\em imbedding} of shrunk random walks 
$B_m(k 2^{-2m})$ into the Wiener process $W(t)$. 

%..............................

\begin{lemma} \label{l-imb}
(a) For any $C \ge 3/2$, $K'>K>0$, and any fixed $m \ge m_0(C,K,K')$ 
there exist random time instants $t_m(k) \in [0, K']$ such that 
\[ \Pb W(t_m(k)) = B_m(k 2^{-2m}),~~~0 \le k 2^{-2m}\le K \rb 
 \ge 1 - 4 (K 2^{2m})^{1-C}, \]
where 
\beq 
\Pb \max _{0 \le k 2^{-2m} \le K} |t_m(k) - k 2^{-2m}| \ge 
 \sqrt {18CKm} \: 2^{-m} \rb \le 4 (K 2^{2m})^{1-C}. \label{eq-tmkk} \eeq 
Moreover, if $\delta $ is such that $0<\delta <1$, $C>2/\delta $, and $m 
\ge m_1(\delta ,C,K,K')$, then we also have  
\beq \Pb \max _{1 \le k 2^{-2m} \le K} |t_m(k) - t_m(k-1) - 2^{-2m}| \ge 
3C 2^{-2m(1-\delta )} \rb 
 \le {K \over 10} 2^{-2m(\delta C -2)} + 4 (K 2^{2m})^{1-C} . 
\label{eq-dtmk} \eeq 

(b)  With probability $1$, for any $K'>K>0$, $0<\delta <1$, and for all but 
finitely many $m$ there exist random time instants $t_m(k) \in [0, K']$ such 
that 
\[ W(t_m(k)) = B_m(k 2^{-2m}) ~~~(0 \le k 2^{-2m} \le K), \]
where 
\[
\max _{0 \le k 2^{-2m} \le K} |t_m(k) - k 2^{-2m}| \le \sqrt {27Km} \: 2^{-m}, 
\]
and 
\[ \max _{1 \le k 2^{-2m} \le K} |t_m(k) - t_m(k-1) - 2^{-2m}| \le 
(7/\delta ) 2^{-2m(1-\delta )}. \]  
\end{lemma}
                                                 
{\sc Proof.}
By Lemma \ref{l-tmn}(a), fixing an $m \ge m_0(C,K,K')$, on a subset $A_m$ of 
the sample space with $\Pb A_m \rb \ge p_m = 1 - 4 (K 2^{2m})^{1-C}$, one has 
\beq \max _{0 \le k 2^{-2m} \le K} |2^{-2n} T_{m,n}(k) - k 2^{-2m}| 
 < \sqrt {18CKm} \: 2^{-m}, \label{eq-tmnk} \eeq  
for each $n>m$. In particular, the time instants $2^{-2n} T_{m,n}(k)$ are
bounded from below by $0$ and from above by $K+\sqrt {18CKm} \: 2^{-m} 
\le K'$. (Assume that $m_0(C,K,K')$ is chosen so.)   

Applying a truncation $t_{m,n}^*(k) = \min \lb K', 2^{-2n} T_{m,n}(k) \rb $, 
for each $k$, $0 \le k 2^{-2m} \le K$, we get a sequence in $n$ bounded over 
the whole sample 
space, equal to the original one for $\omega \in A_m$. It follows from 
the classical Weierstrass theorem \cite[Section 2.42]{Rudin 1977}, that 
every bounded sequence of real numbers contains a convergent subsequence. 
To be definite, let us take the lower limit  
\cite[Section 3.16]{Rudin 1977} of the sequence:
\beq t_m(k) = \liminf _{n \to \infty } t_{m,n}^*(k). \label{eq-tmk} \eeq 
Then $t_m(k) \in [0, K']$. 

By Theorem \ref{thm-w}, with probability $1$ the sample-functions of $B_n(t)$ 
uniformly converge to the corresponding sample-functions of the 
Wiener process that are uniformly continuous on [0, K']. (A continuous 
function on a closed interval is uniformly continuous \cite[Section 4.19]
{Rudin 1977}.) Thus (\ref{eq-bsamec}) implies that for each $k$, $0 \le 
k 2^{-2m} \le K$, we have  $W(t_m(k)) = B_m(k 2^{-2m})$, with probability at 
least $p_m$ (on the set $A_m$ where the truncated sequences coincide with 
the original ones). 

To show it in detail, take any $\epsilon > 0$, any $k$ $(0 \le k 2^{-2m} 
\le K)$, and a subsequence $t_{m,n_i}^*(k)$ converging to $t_m(k)$ as $i \to 
\infty $. Then 
\beqan 
& & |W(t_m(k)) - B_m(k 2^{-2m})| = |W(t_m(k)) - B_{n_i}(2^{-2n_i} 
 T_{m,n_i}(k))| \\
& \le & |W(t_m(k)) - W(2^{-2n_i} T_{m,n_i}(k))| 
 + |W(2^{-2n_i} T_{m,n_i}(k)) - B_{n_i}(2^{-2n_i} T_{m,n_i}(k))| \\
& < & \epsilon /2 + \epsilon /2 = \epsilon ,
\eeqan
where the last inequality holds on the set $A_m$, for all but 
finitely many $n_i$. Since $\epsilon $ was arbitrary, it follows that 
$|W(t_m(k)) - B_m(k 2^{-2m})| = 0$ on $A_m$. 

Further, taking a limit in (\ref{eq-tmnk}) with $n=n_i$ as $i \to \infty $ 
(on the set $A_m$), one obtains (\ref{eq-tmkk}). Also, taking a similar limit 
in Lemma \ref{l-delta}(a), $2^{-2n} \tau _{m,n_i}(k) \to t_m(k)-t_m(k-1)$
on the set $A_m$, and (\ref{eq-dtmk}) follows. 

The statements in (b) can be obtained similarly as in (a), applying Lemmas 
\ref{l-tmn}(b) and \ref{l-delta}(b), or from (a) by the Borel--Cantelli lemma. 
$\Box $

%.............................

\vspace{10pt}

We mention that for any $k \ge 0$ and $m \ge 0$, the sequence $2^{-2n} 
T_{m,n}(k)$ in fact converges to $t_m(k)$ with probability $1$ as $n \to 
\infty $. However, a ``natural'' proof of this fact requires the martingale 
convergence theorem mentioned above, before Lemma \ref{l-delta}, a tool of 
more advanced nature than the ones we use in this paper. 

Next we want to show that the random time instants $s_m(k)$ of the Skorohod 
imbedding (\ref{eq-smk}) and the $t_m(k)$'s defined in (\ref{eq-tmk}) are 
essentially the same. This requires a recollection of some properties of 
random walks. 

We want to estimate the probability that with given positive integers $j$, 
$x$, $u$ and $r$ a random walk $S_i$ goes from a 
point $|S_j|=x$ to $|S_{j+k}|=x+y$ so that $|S_{j+i}| < x+y$ while $1 \le i 
< k$ for some $y \le r$ and $k \ge u$, where $k$, $y$, and $i$ 
are also positive integers. 

The first passage distribution given in \cite[Section III,7]{Feller 1968} can be 
applied here: 
\[ \Pb S_0=0, S_i < y ~(1 \le i < k), S_k=y \rb = {y \over k} {k \choose 
 (k+y)/2} 2^{-k} . \]
Hence, by Theorem \ref{thm-clt1}, 
\[ \Pb |S_j|=x, |S_{j+i}| < x+y ~(1 \le i < k), |S_{j+k}|=x+y ~\mbox{ for 
some }~ y \le r \rb  \]   
\beqan
& \le & \sum _{y=1}^r {y \over k} {k \choose (k+y)/2} 2^{-k} \\
& \le & (1+\epsilon ){r \over k} \left( \Phi (r/\sqrt k) - \Phi (0) \right) \\
& \le & (1+\epsilon ){r \over k} {r \over \sqrt k} {1 \over \sqrt{2 \pi }} 
 \le {r^2 \over k^{3/2}}, 
\eeqan
where $\epsilon > 0$ is arbitrary, say equals $1$, and $k \ge k_0$. 

So the larger the value of $k$ is, the smaller estimate of the probability we 
get. Thus for all positive integers $j$, $x$, $r$, and $u \ge k_0$,   
\beq 
\Pb |S_j|=x, |S_{j+i}| < x+y ~(1 \le i < k), |S_{j+k}|=x+y ~\mbox{ for 
some }~ y \le r,~k \ge u \rb \le r^2/u^{3/2}, \label{eq-fed} \eeq 
independently of the values of $j$ and $x$. 

%..............................

\begin{thm} \label{thm-skor} 
The stopping times $s_m(k)$ $(k \ge 0)$ of the Skorohod imbedding are equal 
to the time instants $t_m(k)$ of the imbedding defined in Lemma \ref{l-imb} 
on the set $A_m$ of the sample space given by 
(\ref{eq-am}), with the possible exception of a zero probability 
subset. 

Therefore all statements in Lemma \ref{l-imb} hold when $s_m(k)$ replaces 
$t_m(k)$. 
\end{thm} 

{\sc Proof.}
Fix an $m \ge m_0(C,K,K')$, where $m_0(C,K,K')$ is the same as in Lemma 
\ref{l-imb}. Let the subset $A_m$ of the sample space be given by 
(\ref{eq-am}). 

Take $k=1$ first. Since $s_m(1)$ is the smallest time instant where $|W(t)|$ 
is equal to $2^{-m}$, and $|W(t_m(1))|=2^{-m}$ on the set $A_m$, it follows 
that $s_m(1) \le t_m(1)$ on $A_m$. We want to show that on $A_m$ the event 
$\lb s_m(1) < t_m(1) \rb $ has zero probability. 
 
Indirectly, let us suppose that $\delta _m = t_m(1)-s_m(1) > 0$ on a subset 
$C_m$ of $A_m$ with positive probability. By (\ref{eq-bsamec}), the first 
time instant where $|B_n(t)|$ equals $|B_m(2^{-2m})|=2^{-m}$ is $2^{-2n} 
T_{m,n}(1)$ $(n>m)$. So $|B_n(t)|<2^{-m}$ if $0 \le t < 2^{-2n}T_{m,n}(1)$. 
On the other hand, by (\ref{eq-asw}), $2^{-m}-n2^{-n/2} \le |B_n(s_m(1))| < 
2^{-m}$ for $n \ge N_1(\omega )$ on a probability $1$ $\omega $-set. 
(Remember that $|W(s_m(1))|=2^{-m}$.) 

Since $\delta_m > 0$ on the set $C_m$, there exists an $N_2(\omega )$ such 
that $n2^{-n/2}<\delta _m /2$ for $n \ge N_2(\omega )$. 

By (\ref{eq-tmk}), $t_m(1)= \liminf _{n \to \infty} 2^{-2n} T_{m,n}(1)$ on 
the set $A_m$. The properties of the lower limit \cite[Section 3.17]{Rudin 
1977} imply that on the subset $C_m$ there exists an $N_3(\omega )$ such that 
$2^{-2n} T_{m,n}(1) > t_m(1) - \delta_m /2$ for $n \ge N_3(\omega )$. 

Set $N(\omega )= \max \lb N_1(\omega ), N_2(\omega ), N_3(\omega ) \rb $ 
for $\omega \in C_m$. 
Since $B_n(t) = 2^{-n} \tilde S_n(t 2^{2n})$, the statements above imply 
that on the set $C_m$ the random walk $\tilde S_n(t)$ have the following 
properties for $n \ge N(\omega )$: 

(a) $| \tilde S_n(s_m(1)2^{2n})| \ge 2^{n-m} - n 2^{n/2}$, 

(b) $| \tilde S_n(t)| < 2^{n-m}$ for $s_m(1)2^{2n} \le t < T_{m,n}(1)$, 
where $T_{m,n}(1) - s_m(1) > (\delta _m /2) 2^{2n} > n 2^{3n/2}$, 

(c) $| \tilde S_n(T_{m,n}(1))| = 2^{n-m}$. 

Let $D_{m,n}$ denote the subset of $C_m$ on which (a), (b), and (c) hold 
for a fixed $n$. Since $D_{m,n} \subset D_{m,n+1}$ for each $n$, by the 
continuity property of probability \cite[Section 11.3]{Rudin 1977}, we have 
$\lim _{n \to \infty} \Pb D_{m,n} \rb = \Pb C_m \rb > 0$. 
This implies that there exists an integer $n_0$ such that 
$\Pb D_{m,n} \rb \ge {1 \over 2} \Pb C_m \rb > 0$ holds for all 
$n \ge n_0$ (say). In other words, for all large enough values of $n$, the 
probability of the event that (a), (b), and (c) hold simultaneously is 
larger than a fixed positive number. 

To get a contradiction, we apply (\ref{eq-fed}) to $\tilde S_n(t)$, with  
$r=n 2^{n/2}$ and $u=n 2^{3n/2}$. Theorem \ref{thm-clt1}, that 
was used to deduce (\ref{eq-fed}), still applies since $r=o(u^{2/3})$, 
i.e. $r/\sqrt u=o(u^{1/6})$. Now the first passage time when $|\tilde S_n(t)|$ 
hits $2^{-2m}$ is $T_{m,n}(1)$. Thus the probability that $\tilde S_n(t)$ 
satisfies (a), (b), and (c) simultaneously is less than or equal to 
\[ {r^2 \over u^{3/2}} = {(n 2^{n/2})^2 \over (n 2^{3n/2})^{3/2}} = 
 {\sqrt n \over 2^{5n/4}}, \] 
which goes to zero as $n \to \infty $. This contradicts the statement above 
that for all large enough value of $n$, the event that (a), (b), and (c) hold 
has a probability larger than a  fixed positive number. This proves 
the lemma for $k=1$: $s_m(1)=t_m(1)$ on the set $A_m$, with the possible 
exception of a zero probability subset. 

For $k>1$, one can proceed by induction. Assume that $s_m(k-1)=t_m(k-1)$ holds 
on $A_m$ except possibly for a subset of probability zero. The proof that 
then $s_m(k)=t_m(k)$ holds as well is essentially the same as the proof of 
the case $k=1$
above. It is true because on one hand $s_m(k)$ is defined recursively in 
(\ref{eq-smk}), using $s_m(k-1)$, the same way as $s_m(1)$ is defined. 
On the other hand, by (\ref{eq-delk}), $T_{m,n}(k)=T_{m,n}(k-1)+
\tau _{m,n}(k)$, where the $\tau _{m,n}(k)$ is defined the same way 
as $\tau _{m,n}(1)=T_{m,n}(1)$. Also, remember that on the set $A_m$, 
$t_m(j)=\liminf _{n \to \infty} T_{m,n}(j)$ for $j=k-1$ or $j=k$. $\Box $

%-------------------------------------------------------------------------

\section{Some Properties of the Wiener Process}

Theorem \ref{thm-w} above indicates that the sample-functions of the Wiener 
process are arbitrarily close to the sample-functions of $B_n(t)$ if $n$ is 
large enough, with probability $1$. The sample-functions of $B_n(t)$ are 
broken lines that have a chance of
$1/2$ to turn and have a corner at any multiple of time $1/2^{2n}$, so at 
more and more instants of time as $n \to \infty $. Moreover, the magnitude 
of the slopes of the line segments that make up the graph of $B_n(t)$ is 
\[ {1/2^n \over 1/2^{2n}} =  2^n \to \infty ~~\mbox{as}~~ n \to  \infty . \]
Therefore one would suspect that the 
sample-functions of the Wiener process are typically nowhere differentiable. 
As we will see below, this is really true. Thus typical sample-functions 
of the Wiener process belong to the ``strange'' class of the everywhere 
continuous but nowhere differentiable functions. 

%.............................

\begin{thm} \label{thm-ndiff}
With probability $1$, the sample-functions of the Wiener process are nowhere 
differentiable. 
\end{thm}

{\sc Proof.}
It suffices to show that with probability $1$, the sample-functions are 
nowhere differentiable on any interval $[0,K]$. Put $K_0=(3/2)K >0$ (say).   
Then with probability $1$, for all sample-functions and 
for all but finitely many $m$ there exist time instants $t_m(k)$ 
$(0 \le k 2^{-2m} \le K_0)$ with the properties described in Lemma 
\ref{l-imb}(b). In particular, 
\[ \max _{0 \le k 2^{-2m} \le K_0} t_m(k) \ge K_0 - \sqrt { 27 K_0 m} 
\: 2^{-m} > K \] 
if $m$ is large enough. 

Fix an $\omega $ in this probability $1$ subset of the sample 
space. This defines a specific sample-function of $W(t)$ and specific values 
of the random time instants $t_m(k)$. 
(To simplify the notation, in this proof we suppress the argument $\omega $.)
Then choosing an arbitrary point $t \in [0,K]$, for each $m$ large enough, 
one has $t_m(k-1) \le t < t_m(k)$ for some $k$, $0 < k 2^{-2m} \le K_0$. 
Taking for instance $\delta = 1/4$ in Lemma \ref{l-imb}(b), we get $t_m(k) 
- t_m(k-1) \le 29 \cdot 2^{-(3/2)m}$ and 
\[ |W(t_m(k)) - W(t_m(k-1))| = |B_m(k 2^{-2m}) - B_m((k-1) 2^{-2m})| = 2^{-m}. 
\]
 
Set $t_m^* = t_m(k)$ if $|W(t) - W(t_m(k))| \ge |W(t) - W(t_m(k-1))|$ and 
$t_m^* = t_m(k-1)$ otherwise. Then $|W(t) - W(t_m^*)| \ge (1/2) 2^{-m}$. So  
$|t_m^* - t| \le 29 \cdot 2^{-(3/2)m} \to 0$ and  
\[ \left| {W(t_m^*) - W(t) \over t_m^* - t} \right| \ge {(1/2) 2^{-m} \over 
 29 \cdot 2^{-(3/2)m}} = {1 \over 58} 2^{m/2}  \to \infty , \]
as $m \to \infty $. This shows that the given sample-function cannot be 
differentiable at any point $t \in [0,K]$. $\Box $

%.............................

\vspace{10pt}

It has important consequences in the definition of stochastic integrals that, 
as shown below, 
the graph of a typical sample-function of the Wiener process has infinite 
length. In general, (the graph of) a function $f$ defined on an interval 
$[a,b]$ has {\em finite length} (or $f$ is said to be of {\em bounded 
variation} on $[a,b]$) if there exists a finite constant $c$ such that for 
any {\em partition} $a=x_0<x_1< \cdots <x_{n-1}<x_n=b$, the sum of the 
absolute values of the corresponding changes does not exceed $c$: 
\[ \sum _{j=1}^n |f(x_j) - f(x_{j-1})| \le c . \] 
The smallest $c$ with this property is called the {\em total variation} of 
$f$ over $[a, b]$, denoted $V(f(t), ~ a \le t \le b)$. 
Otherwise we say that the graph has {\em infinite length}, or $f$ is of 
{\em unbounded variation} on $[a,b]$. 

First let us calculate the total variation of a sample-function of $B_m(t)$ 
over an interval $[0,K]$. Each sample-function of $B_m(t)$ over $[0,K]$ is a 
broken line that consists of $K 2^{2m}$ line segments with changes of 
magnitude $2^{-m}$. So for any sample-function of $B_m(t)$, 
\beq  V(B_m(t),~0 \le t \le K) = K 2^{2m} 2^{-m} = K 2^m , \label{eq-lbm} \eeq
which tends to infinity as $m \to \infty $. 

\begin{lemma} \label{l-ubv} 
For any $K'>0$, the sample-functions of the Wiener process over $[0, K']$ 
have infinite length (i.e. are of unbounded variation) with probability $1$. 
\end{lemma}

{\sc Proof.} 
By Lemma \ref{l-imb}, for any $C \ge 3/2$, $K'>K>0$, and $m \ge m_0(C,K,K')$ 
there exist time instants $t_m(k) \in [0, K']$ such that 
\beq \Pb W(t_m(k)) = B_m(k 2^{-2m}),~~~0 \le k 2^{-2m}\le K \rb 
 \ge 1 - 4 (K 2^{2m})^{1-C}. \label{eq-imb1} \eeq

For each $m \ge 0$ define the following event: 
\[ C_m = \lb V(W(t),~0 \le t \le K')  < K 2^m \rb . \] 
Then $C_m \subset C_{m+1}$ for any $m \ge 0$. 

For any sample-function of 
$W(t)$, take the partition $0=t_m(0)<t_m(1)< \cdots <t_m(K 2^{2m})$. (To 
alleviate the notation, we suppress the dependence on $\omega $.) By 
(\ref{eq-imb1}), for any $m \ge m_0(C,K,K')$, the sum of the corresponding 
absolute changes is equal to $K 2^{2m} 2^{-m}=K 2^m$, with probability at 
least $1-4(K 2^{2m})^{1-C}$. 

This shows that then $\Pb C_m \rb < 4 (K 2^{2m})^{1-C}$. Take the event 
\[ C_{\infty } = \lb V(W(t),~0 \le t \le K') < \infty \rb . \]
The continuity property of probability implies that $\Pb C_m \rb \to 
\Pb C_{\infty } \rb $ as $m \to \infty $, that is, $\Pb C_{\infty } \rb = 0$. 
$\Box $ 

%.............................

\vspace{10pt}

The next lemma shows a certain uniform continuity property of the Wiener 
process. An interesting consequence of the lemma is that for any $u>0$ the 
probability that $|W(t)-W(s)| \ge u$ holds for some $s,t \in [0,K]$, $|t-s| 
\le h$ can be made arbitrarily small if a small enough $h$ is chosen. 
More accurately, the lemma shows that only with small probability can the 
increment of the Wiener process be larger than $c\sqrt h$ if the constant 
$c$ is large enough. Now $\sqrt h$ is much larger than $h$ for small 
values of $h$, so this also indicates why sample-functions of the 
Wiener process are not differentiable. At the same time it gives a rough 
measure of the so-called {\em modulus of continuity} of the process.  
Basically, the proof relies on Theorem \ref{thm-clt1}a and Theorem \ref{thm-w}. 

\begin{lemma}  \label{l-uconp} 
For any $K>0$, $0<\delta <1$, and $u>0$ there exists an $h_0(K, \delta , u)>0$ 
such that 
\beq \Pb \max _{s,t \in [0,K],~|t-s| \le h} |W(t) - W(s)| \ge u \rb \le 
 7 e^{-{u^2 \over 2h} (1-\delta )} , \label{eq-uconp} \eeq
for all positive $h \le h_0(K, \delta , u)$. 
\end{lemma}

{\sc Proof.}
First we choose a large enough $C \ge 3/2$ 
such that $2/(C-1) < \delta /2$. For instance, $C=1+(6/\delta )$  
will do. 

By (\ref{eq-pw}), the probability in (\ref{eq-uconp}) cannot exceed 
\beq 6(K 2^{2n})^{-6/\delta } + \Pb \max _{0 \le s \le K-h} \max _{s \le t 
 \le s+h} |B_n(t) - B_n(s)| \ge u - 2n 2^{-n/2} \rb , \label{eq-dbn} \eeq
for $n \ge n_0(K,\delta )$. (Remember that $1-C=-6/\delta $ now.) 

By definition, $B_n(t)=2^{-n} \tilde S_n(t 2^{2n})$ for $t \ge 0$. For each 
$s \le t$ from $[0, K]$ and $n \ge n_0(K,\delta )$ take the integers 
$s_n=\lceil s 2^{2n} \rceil $ 
and $t_n=\max \lb s_n, \lfloor t 2^{2n} \rfloor \rb $. ($\lceil x \rceil $ 
denotes the smallest integer greater than or equal to $x$, while $\lfloor x 
\rfloor $ is the largest integer smaller than or equal to $x$.) 

Then $|t_n - t 2^{2n}| \le 1$ and so $|\tilde S_n(t_n) - \tilde S_n(t 2^{2n})| 
\le 1$, similarly for $s_n$. Moreover, $0 \le t_n - s_n \le h 2^{2n}$ 
if $0 \le t-s \le h$. Hence (\ref{eq-dbn}) does not exceed 
\beq 6(K 2^{2n})^{-6/\delta } + \Pb \max _{0 \le j \le K 2^{2n}} \max _{0 \le 
 k \le h 2^{2n}} |\tilde S_n(j+k) - \tilde S_n(j)| \ge 2^n (u - 2n 2^{-n/2}) 
 - 2 \rb , \label{eq-dsn} \eeq
for $n \ge n_0(K,\delta )$. 

The distribution of $\tilde S_n(j+k) - \tilde S_n(j)$ above is the same as
the distribution of a random walk $S(k)$, for any value of $k \ge 0$, 
independently of $j \ge 0$. Also, the largest possible value of $|S(k)|$ is 
$k$. Therefore by Theorem \ref{thm-clt1}a, the inequality (\ref{eq-asy}), and 
the crude estimate (\ref{eq-crude}),
\beqan \lefteqn 
{\Pb \max _{0 \le k \le h 2^{2n}} |\tilde S_n(j+k) - \tilde S_n(j)| 
\ge 2^n (u - 2n 2^{-n/2} - 2 \cdot 2^{-n}) \rb } \\
& \le & \Pb \max _{u \sqrt {1-\delta /2} \: 2^n \le k \le h 2^{2n}} {|S(k)| 
 \over \sqrt k} \ge {u \over \sqrt h} \sqrt {1-\delta /2} \rb 
 \le h 2^{2n} e^{-{u^2 \over 2h} (1-\delta /2)} .
\eeqan
Here it was assumed that $2n 2^{-n/2} + 2\cdot 2^{-n} \le u (1 - \sqrt {1-
\delta /2})$, which certainly holds if $n \ge n_1(K,\delta ,u) \ge n_0(K,
\delta )$. Also, we assumed that ${u \over \sqrt h} \sqrt {1-\delta /2} 
\ge 3/\sqrt {2 \pi}$, see (\ref{eq-ldt1}), which is 
true if $h$ is small enough, depending on $\delta $ and $u$. 

Consequently, applying the crude estimate (\ref{eq-crude}) again for 
(\ref{eq-dsn}), we obtain 
\beqan  \lefteqn
{\Pb \max _{s,t \in [0,K],~|t-s| \le h} |W(t) - W(s)| \ge u \rb } \\ 
& \le & 6(K 2^{2n})^{-6/\delta } + K 2^{2n} h 2^{2n} e^{-{u^2 \over 2h} 
 (1-\delta /2)} \\ 
& = & 6 e^{-{6 \over \delta } (\log K + 2n \log 2)} + Kh e^{4n \log 2 - 
 {u^2 \over 2h} (1-\delta /2)}. 
\eeqan 

Now we select an integer $n \ge n_1(K,\delta , u)$ such that 
$-{6 \over \delta } (\log K + 2n \log 2) \le -{u^2 \over 2h}$.  
The choice 
\[ n = \left\lceil {1 \over 2 \log 2} \left( {u^2 \over 2h} {\delta \over 6} 
 - \log K \right) \right\rceil  \] 
will do if $h$ is small enough, $0 < h \le h_0(K,\delta , u)$, so that 
$n \ge n_1(K,\delta , u) \ge 2$. Then $n \le {3 \over 2} {1 \over 2 \log 2} 
\left( {u^2 \over 2h} {\delta \over 6} - \log K \right)$ holds as well. 

With this $n$ we have $ 4n \log 2 \le {u^2 \over 2h} \delta /2 + \log 
(K^{-3})$, and so 
\beqan  \lefteqn
{\Pb \max _{s,t \in [0,K],~|t-s| \le h} |W(t) - W(s)| \ge u \rb } \\ 
& \le & 6 e^{-{u^2 \over 2h}} + Kh K^{-3} e^{- {u^2 \over 2h}(1-
 \delta /2 - \delta /2)} \\  
& \le & (6 + h/K^2) e^{-{u^2 \over 2h} (1-\delta )} . 
\eeqan 

If $K \ge 1$, then $h/K^2 \le 1$ and (\ref{eq-uconp}) follows. If $K<1$, the 
maximum in (\ref{eq-uconp}) cannot exceed the maximum over the interval 
$[0, 1]$. Then taking $h_0(K,\delta ,u)=h_0(1,\delta ,u)$, (\ref{eq-uconp}) 
follows again. $\Box $ 

%-------------------------------------------------------------------------

\section{A Preview of Stochastic Integrals}

To show how stochastic integrals come as natural tools when working with 
differential equations including random effects, and 
what kind of problems arise when one wants to define them, let us start 
with the simplest ordinary differential equation 
\[ x'(t) = f(t)~~~  (t \ge 0), \] 
where $f$ is a continuous function. If $x(0)$ is given, its unique solution 
can be obtained by integration,  
\[ x(t) - x(0) = \int_0^t f(s)~ds~~~(t \ge 0) . \]

Now we modify this simple model by introducing a random term, very customary  
in several applications: 
\[ x'(t) = f(t) + g(t) W'(t)~~~(t \ge 0), \]
where $f$ and $g$ are continuous random functions and $W'(t)$ is the so-called 
{\em white noise} process. Now we know from Theorem \ref{thm-ndiff} that 
$W'(t)$ does not exist (at least not in the ordinary sense), but after 
integration we may get some meaningful solution, 
\[ x(t) - x(0) = \int _0^t f(s)~ds + \int _0^t g(s)~dW(s)~~~(t \ge 0) . \] 
The second integral here is what one wants to call a stochastic integral 
if it can be defined properly.  
 
A natural idea to define such a stochastic integral is to define it as a 
{\em Riemann--Stieltjes integral} \cite[Chapter 6]{Rudin 1977} for each 
sample-function separately. It means that one takes partitions 
$0=s_0<s_1<\cdots <s_{n-1}<s_n=t$, and Riemann--Stieltjes sums 
\[ \sum _{k=1}^n g(u_k) (W(s_k) - W(s_{k-1})) , \]
where $u_k \in [s_{k-1}, s_k]$ is arbitrary. (We suppress the argument 
$\omega $ that would refer to a specific sample-function in order to 
alleviate the notation.) Then one would hope that as the norm of the 
partition $\| {\cal P} \| = \max _{1 \le k \le n} |s_k - s_{k-1}|$ tends 
to $0$, the Riemann--Stieltjes sums converge to the same limit when 
fixing a specific point $\omega $ in the sample space. 

One problem is that it cannot happen to all continuous random functions 
$g$. The reason is that $W(s)$ has unbounded variation over the interval 
$[0,t]$ ---as we saw it in Lemma \ref{l-ubv}. The random function $g$ could 
be chosen so that a Riemann--Stieltjes sum gets arbitrary close to the 
total variation, which is $\infty $. 
Naturally, this is the case with not only the Wiener process, but with 
any process whose sample functions have unbounded variation, see e.g. 
\cite[Section I.7]{Protter 1992}.  

But there is another problem connected to the choice of the points $u_k 
\in [s_{k-1}, s_k]$ in the Riemann--Stieltjes sums above. This choice 
unfortunately does matter, not like in the case of ordinary integration. 
The reason is again the unbounded variation of the sample-functions. 
The easiest way 
to illustrate it is using {\em discrete stochastic integrals}, that 
is, sums of random variables. (Such a sum is essentially the same as a 
Riemann--Stieltjes sum above.)  

So let $S_0=0$, $S_n=\sum _{k=1}^n X_k$ is a (simple, symmetric) random 
walk, just like in Section 1. In the following examples $S_n$ will play 
the role of the function $g(t)$ above, and the white noise process $W'(t)$ 
is substituted by the increments $X_n$. In the first case (that 
corresponds to an {\em It\^o-type stochastic integral\/}), we define the 
discrete stochastic integral as $\sum _{k=1}^n S_{k-1} X_k$. Observe that 
in this case the integrand is always taken at the left endpoint of the 
subintervals. A usual reasoning behind this is that $X_k$ gives the ``new 
information'' in each term, while the integrand $S_{k-1}$ depends only on 
the past, that is, {\em non-anticipating\/}: independent of the 
future values $X_k, X_{k+1}, \ldots $.   

This discrete stochastic integral can be evaluated explicitly as 
\beqan 
\sum _{k=1}^n S_{k-1} X_k & = & \sum _{k=1}^n S_{k-1} (S_k - S_{k-1}) \\
& = & {1 \over 2} \sum _{k=1}^n (S_k^2 - S_{k-1}^2) - {1 \over 2} 
\sum _{k=1}^n (S_k - S_{k-1})^2 = {S_n^2 \over 2} - {n \over 2} .
\eeqan
Here we used that the first resulting sum telescopes and $S_0^2=0$, while 
each term  $(S_k-S_{k-1})^2$ in the second resulting sum is equal to $1$. 
The interesting feature of the result is that it contains the 
non-classical term $-n/2$. The ``non-classical'' phrase refers to the fact 
that $\int _0^{s_n} s~ds = {s_n}^2/2$. Altogether, this formula is a special 
case of the important {\em It\^o formula}, one of our main subjects from 
now on. 

Of course, it is also interesting to see what happens if the integrand is 
always evaluated at the right endpoints of the subintervals: 
\beqan 
\sum _{k=1}^n S_k X_k & = & \sum _{k=1}^n S_k (S_k - S_{k-1}) \\
& = & {1 \over 2} \sum _{k=1}^n (S_k^2 - S_{k-1}^2) + {1 \over 2} 
\sum _{k=1}^n (S_k - S_{k-1})^2 = {S_n^2 \over 2} + {n \over 2} .
\eeqan
Note that the non-classical term is $+n/2$ here. 

Taking the arithmetical average of the two formulae above we obtain a 
{\em Stratonovich-type stochastic integral}, which does not contain a 
non-classical term: 
\[ \sum _{k=1}^n {S_{k-1}+S_k \over 2} X_k = \sum _{k=1}^n S( k-
{1 \over 2}) X_k = {S_n^2 \over 2} . \]
On the other hand, this type of integral has other disadvantages compared 
to the It\^o-type one, resulting from the fact that here the 
integrand is ``anticipating'', not independent of the future. 

After showing these (and other) examples in a seminar, P. R\'ev\'esz 
asked the question if there is a general method to evaluate discrete 
stochastic integrals of the type $\sum_{k=1}^n f(S_{k-1}) X_k$ in closed 
form, where $f$ is a given function defined on the set of 
integers \ZZ . In other words, does there exist a discrete It\^o formula 
in general? The answer is yes, and fortunately it is quite elementary to 
see.

But before turning to this, let us see the relationship of such a 
formula to an alternative way of defining certain stochastic integrals. 
This important type of stochastic integrals is $\int _0^K f(W(s))~dW(s)$, 
where $K>0$ and $f$ is a continuously differentiable function. In other 
words, the integrand is a smooth function of the Wiener process. The 
traditional definition of the It\^o-type integral in this case goes 
quite similarly to the Riemann--Stieltjes integral. 

Take an arbitrary partition ${\cal P} = \lb 0=s_0, s_1, 
\ldots ,s_{n-1}, s_n=K \rb $ on the time axis, and 
a corresponding Riemann-Stieltjes sum, evaluating the function always at 
the left endpoints of the subintervals, 
\[ \sum _{k=1}^n f(W(s_{k-1})) ~(W(s_k) - W(s_{k-1})) . \]
This sum is a random variable, corresponding to the given partition. It 
can be proved that these random variables converge e.g. {\em in probability} 
to a certain random variable $I$, as the norm of the partition goes to $0$. 
This random variable $I$ is then called the It\^o integral. We mention 
that ``in probability'' convergence means that for any $\epsilon >0$ there 
exists a $\delta >0$ such that  
\[ \Pb |I - \sum _{k=1}^n f(W(s_{k-1})) ~(W(s_k) - W(s_{k-1}))| \ge \epsilon 
\rb  < \epsilon , \] 
as $\| {\cal P} \| < \delta $. 

The alternative method that we will follow in this paper is better suited to 
the relationship between the Wiener process and random walks discussed above. 
Mathematically, it somewhat reminds a {\em Lebesgue--Stieltjes integral} 
\cite[Chapter 11]{Rudin 1977}. 
The idea is that we first take a dyadic partition on the spatial axis, 
each subinterval having the length $2^{-m}$, where $m$ is a non-negative 
integer. Then we determine the corresponding first passage times $s_m(1), 
s_m(2), \ldots $ of the Skorohod imbedding as explained above. These time 
instants can be considered as a random partition on the time axis that in 
general depends on the considered sample-function.   

By Lemma \ref{l-imb}b and Theorem \ref{thm-skor}, with probability $1$, for 
any $K'>0$ and for all but finitely many $m$, each $s_m(k)$ lies in the 
interval $[0, K']$ and $W(s_m(k)) = B_m(k 2^{-2m})$, $0 \le k 2^{-2m} \le 
K$. The shrunk random walk $B_m(t)$ can be expressed in terms of ordinary 
random walks by (\ref{eq-bmt}) as $B_m(k 2^{-2m}) = 2^{-m} \tilde S_m(k)$. 
Now our definition of the It\^o integral will be 
\beq \lim _{m \to \infty } \sum _{k=1}^{K 2^{2m}} f(W(s_m(k-1)))~(W(s_m(k)) 
- W(s_m(k-1))) . \label{eq-idef1} \eeq

We will show later that this sum, which can be evaluated for each 
sample-function separately, converges with probability $1$. Our method 
will be to find an other form of this sum by a discrete It\^o formula and 
to apply the limit to the equivalent form so obtained. 

%-------------------------------------------------------------------------

\section{A Discrete It\^o Formula}

Let $f$ be a function defined on the set of integers $\ZZ $. First we 
define {\em trapezoidal sums} of $f$ by  
\beq T_{j=0}^k ~f(j) = \epsilon _k \lb {1\over 2} f(0) + \sum _{j=1}^{|k|-1} 
f(\epsilon _k j) ~+ {1\over 2} f(k) \rb , \label{eq-trap} \eeq 
where $k \in \ZZ $ (so $k$ can be negative as well!) and  
\beq \epsilon _k = \left\{ \begin{array}{rl}
1 & \mbox{ if $~k>0$} \\
0 & \mbox{ if $~k=0$} \\
-1 & \mbox{ if $~k<0$ .} 
\end{array} \right. \label{eq-eps} \eeq 
The reason behind the $-1$ factor when $k<0$ is the analogy with 
integration: when the upper limit of the integration is less than the 
lower limit, one can exchange them upon multiplying the integral by $-1$. 

The next statement that we will call a {\em discrete It\^o formula} is a 
purely algebraic one. It is reflected by the fact that though we will 
apply it exclusively for random walks, the lemma holds for any numerical  
sequence $X_r=\pm 1$, irrespective of any probability assigned to them.  

%.............................

\begin{lemma}  \label{l-dito}
Take any function $f$ defined on $\ZZ $, any sequence $X_r=\pm 1$ ($r \ge 
1$), and let $S_0=0$, $S_n=X_1+X_2+\cdots +X_n$ $(n \ge 1)$. Then the 
following statements hold: 

{\sc discrete It\^o formula}
\[ T_{j=0}^{S_n}~f(j) = \sum _{r=1}^n f(S_{r-1}) 
X_r + {1\over 2} \sum _{r=1}^n {f(S_r) - f(S_{r-1}) \over X_r} ,\] 
and 

{\sc discrete Stratonovich formula} 
\[ T_{j=0}^{S_n}~f(j) = \sum _{r=1}^n {f(S_{r-1}) + f(S_r) \over 2} X_r .\]
\end{lemma}

{\sc Proof.} By the definition of a trapezoidal sum, 
\beq T_{j=0}^{S_r}~f(j) - T_{j=0}^{S_{r-1}}~f(j) = X_r {f(S_{r-1}) + f(S_r) 
\over 2}, \label{eq-dito1} \eeq 
since if $S_r-S_{r-1}=X_r$ equals $1$, one has to add a term $(f(S_{r-1}) + 
f(S_r))/2$, while if $X_r=-1$, one has to subtract this term. 

Since $X_r=\pm 1$, the right hand side of (\ref{eq-dito1}) can be written 
as 
\beq T_{j=0}^{S_r}~f(j) - T_{j=0}^{S_{r-1}}~f(j) = f(S_{r-1}) X_r + 
{1\over 2} {f(S_r) - f(S_{r-1}) \over X_r}. \label{eq-dito2} \eeq 

By summing up (\ref{eq-dito2}), respectively (\ref{eq-dito1}), for 
$r=1,2,\ldots ,n$ we obtain the statements of the lemma, since the sum 
telescopes and $T_{j=0}^{S_0}~f(j) = 0$:
\[ \sum _{r=1}^n \left( T_{j=0}^{S_r}~f(j) - T_{j=0}^{S_{r-1}}~f(j) 
\right) = T_{j=0}^{S_n}~f(j) . \] 
$\Box $ 

%.............................

\vspace{10pt}

We need a version of Lemma \ref{l-dito} that can be applied for shrunk 
random walks $B_m(t)$ as well. Therefore we define trapezoidal sums of a 
function $f$  over an equidistant partition with points $x=j \Delta x$, where 
$\Delta x >0$ and $j$ changes over the set of integers $\ZZ $. Here the 
function $f$ is assumed to be defined on the set of real numbers $\RR $. 
So a corresponding trapezoidal sum is 
\beq T_{x=0}^a~f(x)~\Delta x = \epsilon _a \Delta x \lb {1\over 2} 
f(0) + \sum _{j=1}^{(|a|/\Delta x)-1} f(\epsilon _a j \Delta x) ~+ {1\over 2} 
f(a) \rb , \label{eq-trapm} \eeq 
where $a$ is assumed to be an integer multiple of $\Delta x$ and 
$\epsilon_a$ is defined according to (\ref{eq-eps}).  
In the sequel this definition will be applied with $\Delta x = 2^{-m}$. 
We write the corresponding version of Lemma \ref{l-dito} directly for 
shrunk random walks $B_m(t)$, though this lemma is of purely algebraic 
nature as well.  

%.............................

\begin{lemma}  \label{l-ditom}
Take any function $f$ defined on $\RR $, any real $K>0$, and fix a non-negative 
integer $m$. Consider shrunk random walks $B_m(r 2^{-2m})=2^{-m} \tilde 
S_m(r)$ ($r \ge 0$). Then the following statements hold ($\Delta 
x=2^{-m}$, $\Delta t=2^{-2m}$): 

{\sc It\^o case}
\beqa  
T_{x=0}^{B_m(K_m)}~f(x)~\Delta x 
& = & \sum _{r=1}^{\lfloor K/\Delta t\rfloor } 
f(B_m((r-1) \Delta t))~(B_m(r \Delta t) - B_m((r-1) \Delta t)) 
\nonumber \\  
& + & {1\over 2} \sum _{r=1}^{\lfloor K/\Delta t\rfloor } {f(B_m(r \Delta t)) 
- f(B_m((r-1) \Delta t)) \over B_m(r \Delta t) - B_m((r-1) \Delta t) } 
\Delta t ,
\label{eq-ditom} \eeqa 
and 

{\sc Stratonovich case}
\beqa \lefteqn 
{T_{x=0}^{B_m(K_m)}~f(x)~\Delta x} \nonumber \\
& = & \sum _{r=1}^{\lfloor K/\Delta t\rfloor } {f(B_m((r-1) \Delta t)) + 
f(B_m(r \Delta t)) \over 2} (B_m(r \Delta t) - B_m((r-1) \Delta t)) , 
\label{eq-dstratm} \eeqa 
where $K_m=\lfloor K/\Delta t\rfloor \Delta t$. 
\end{lemma}

{\sc Proof.} 
The proof is essentially the same as in case of Lemma \ref{l-dito}, therefore 
omitted. $\Box $

%.............................

\vspace{10pt}

Now recall Lemma \ref{l-imb}b and Theorem \ref{thm-skor}. With 
probability $1$, for any $K'>K$ and for all but finitely many $m$ there 
exist random time instants $s_m(r) \in [0, K']$ (the first passage times 
of the Skorohod imbedding) such that $W(s_m(r))=B_m(r \Delta t)$ and 
\beq \max _{0 \le r \Delta t \le K} |s_m(r) - r \Delta t | \le \sqrt {27Km} 
\: 2^{-m}, \label{eq-edist} \eeq 
going to $0$ as $m \to \infty $. 

In this light the shrunk random walks 
$B_m(t)$ can be replaced by the Wiener process in (\ref{eq-ditom}) and 
(\ref{eq-dstratm}). Then the first sum on the right hand side of 
(\ref{eq-ditom}) becomes exactly the one whose limit as $m \to \infty $ 
is going to be our definition of It\^o integral by (\ref{eq-idef1}). 
Similarly, the right hand side of (\ref{eq-dstratm}) is the one whose 
limit will be our definition of the Stratonovich integral. 

The most important feature of Lemma \ref{l-ditom} is that these limits 
can be evaluated in terms of limits of other, simpler sums. An other gain 
is that after performing the limits, we will immediately obtain the 
important It\^o and Stratonovich formulae for the corresponding types of 
stochastic integrals.

%-------------------------------------------------------------------------

\section{Stochastic Integrals and the It\^o formula}

\begin{thm} Let $f$ be a continuously differentiable function on the set 
of real numbers $\RR $, and $K>0$. For $m \ge 0$ and $k \ge 0$ take the 
first passage times $s_m(k)$ of the Skorohod imbedding of shrunk random 
walks into the Wiener process as defined by (\ref{eq-smk}). Then the sums 
below converge with probability $1$: 

{\sc It\^o integral}
\beq \int _0^K f(W(s))~dW(s) = \lim _{m \to \infty } \sum _{r=1}^{K 2^{2m}} 
f(W(s_m(r-1)))~(W(s_m(r)) - W(s_m(r-1))) , \label{eq-itod} \eeq
and

{\sc Stratonovich integral}
\beqa \lefteqn 
{\int _0^K f(W(s))~\circ ~dW(s) } \nonumber \\ 
& = & \lim _{m \to \infty } \sum _{r=1}
^{K 2^{2m}} {f(W(s_m(r-1))) + f(W(s_m(r))) \over 2}~(W(s_m(r)) - 
W(s_m(r-1))) . 
\label{eq-strd} \eeqa

For the corresponding stochastic integrals we have the following formulae 
as well: 

{\sc It\^o formula}
\beq \int _0^{W(K)} f(x)~dx = \int _0^K f(W(s))~dW(s) + {1\over 2} \int 
_0^K f'(W(s))~ds , \label{eq-itof} \eeq 
and 

{\sc Stratonovich formula}
\beq \int _0^{W(K)} f(x)~dx = \int _0^K f(W(s))~\circ ~dW(s) . 
\label{eq-straf} \eeq 
\end{thm}

{\sc Proof.}
By the It\^o case of Lemma \ref{l-ditom} and the comments made after lemma, 
with probability $1$, for all but finitely many $m$, we have the next 
equation for the sum in (\ref{eq-itod}):  
\beqa \lefteqn 
{\sum _{r=1}^{\lfloor K/\Delta t \rfloor } {f(W(s_m(r-1))) + f(W(s_m(r))) 
\over 2}~(W(s_m(r)) - W(s_m(r-1)))} \nonumber \\ 
& = & T_{x=0}^{W(s_m(\lfloor K/\Delta t \rfloor ))}~f(x)~\Delta x - {1\over 2} 
\sum _{r=1}^{\lfloor K/\Delta t \rfloor } {f(W(s_m(r))) - f(W(s_m(r-1))) 
\over W(s_m(r)) - W(s_m(r-1))} \Delta t , \label{eq-iteq} \eeqa  
where $\Delta x = 2^{-m}$ and $\Delta t = 2^{-2m}$. 

For $t \in [0, K]$ set $t_m=\lfloor t/\Delta t \rfloor \Delta t$. Then 
$|t-t_m| \le \Delta t=2^{-2m}$. By (\ref{eq-edist}), $|t_m - s_m(\lfloor 
t/\Delta t \rfloor)| \le \sqrt {27Km}\: 2^{-m}$ with probability 
$1$ if $m$ is large enough. This implies that 
\beq \max _{0 \le t \le K} |t - s_m(\lfloor t/\Delta t \rfloor)| \to 0 
\label{eq-tlim} \eeq  
with probability $1$ as $m \to \infty $. Further, the sample functions 
of the Wiener process being uniformly continuous on $[0, K']$ with 
probability $1$, one gets that then 
\beq \max _{0 \le t \le K} |W(t) - W(s_m(\lfloor t/\Delta t \rfloor))| \to 0 
\label{eq-xlim} \eeq 
as well.

Particularly, it follows that $W(s_m(\lfloor K/\Delta t \rfloor)) 
\to W(K)$ with probability $1$ as $m \to \infty $. 
On the other hand, the trapezoidal sum $T_{x=0}^a~f(x)~\Delta x$ of a 
continuous function $f$ is a Riemann sum corresponding to the partition 
$\lb 0, {1 \over 2} \Delta x, {3 \over 2} \Delta x, \ldots , a - {3 \over 2} 
\Delta x, a - {1 \over 2} \Delta x, a \rb $. Therefore the trapezoidal sums 
converge to $\int _{x=0}^a f(x)~dx$ as $\Delta x \to 0$. 
These show that for any $\epsilon >0$, 
\beqan \lefteqn 
{\left| \int _0^{W(K)} f(x)~dx - T_{x=0}^{W(s_m(\lfloor K/\Delta t 
\rfloor))}~f(x)~\Delta x \right|} \\
& \le & \left| \int _0^{W(s_m(\lfloor K/\Delta t \rfloor))} f(x)~dx - 
T_{x=0}^{W(s_m(\lfloor K/\Delta t \rfloor))}~f(x)~\Delta x \right| + 
\left| \int _{W(s_m(\lfloor K/\Delta t \rfloor))}^{W(K)} f(x)~dx \right| \\
& < & \epsilon /2 + \epsilon /2 = \epsilon 
\eeqan 
with probability $1$ if $m$ is large enough. That is, the trapezoidal sum 
in (\ref{eq-iteq}) tends to the corresponding integral with probability $1$: 
\beq \lim _{m \to \infty } T_{x=0}^{W(s_m(\lfloor K/\Delta t \rfloor))}~f(x)
~\Delta x = \int _0^{W(K)} f(x)~dx . \label{eq-traco} \eeq 

Now let us turn to the second sum in (\ref{eq-iteq}). By the definition 
of the first passage times, $W(s_m(r)) - W(s_m(r-1)) = \pm 2^{-m} = \pm 
\Delta x$, which tends to $0$ as $m \to \infty $. Hence 
\beq {f(W(s_m(r))) - f(W(s_m(r-1))) \over W(s_m(r)) - W(s_m(r-1))}
= {f(W(s_m(r)) \mp \Delta x) - f(W(s_m(r))) \over \mp \Delta x} . 
\label{eq-difq} \eeq 
We want to show that this difference quotient gets arbitrarily close to
$f'(W(r \Delta t))$ if $m$ is large enough. 

To this end, let us consider the following problem from calculus. If $f$ 
is a continuously differentiable function, $x_m \to x$ and $\Delta x_m 
\to 0$ as $m \to \infty $, let us consider the difference of $f'(x)$ and 
$(f(x_m+\Delta x_m)-f(x_m))/\Delta x_m$. By the mean value theorem, the 
latter diference quotient is equal to $f'(u_m)$, where $u_m$ lies 
between $x_m \to x$ and $x_m+\Delta x_m \to x$. Since $f'$ is continuous, 
this implies that 
\beq {f(x_m+\Delta x_m)-f(x_m)) \over \Delta x_m} \to f'(x) , 
\label{eq-udif} \eeq 
as $m \to \infty $. 

In our present context, $x=W(t)$ and $x_m=W(s_m(\lfloor t/\Delta t 
\rfloor))$, where $0 \le t \le K$. Since the sample functions of $W(t)$ 
are continuous with probability $1$, it follows from the max-min theorem 
that their ranges are contained in bounded intervals. Over such a bounded 
interval the function $f'$ is uniformly continuous, therefore 
(\ref{eq-tlim}), (\ref{eq-xlim}), and (\ref{eq-udif}) imply 
\beq \max _{0 \le t \le K} \left| f'(W(t)) - {f(W(s_m(\lfloor t/\Delta t 
\rfloor)) \mp \Delta x) - f(W(s_m(\lfloor t/\Delta t \rfloor))) \over 
\mp \Delta x} \right| \to 0 \label{eq-wdlim} \eeq 
with probability $1$ as $m \to \infty$. (Remember that now $\Delta x = 
2^{-m}$ and $\Delta t = 2^{-2m}$.) 

Particularly, for any $\epsilon >0$, we have   
\beqa 
& & \max _{1 \le r \Delta t \le K} \left| {f(W(s_m(r))) - f(W(s_m(r-1))) 
\over W(s_m(r)) - W(s_m(r-1))} - f'(W(r \Delta t)) \right| \nonumber \\
& = & \max _{1 \le r \Delta t \le K} \left| {f(W(s_m(r)) \mp \Delta x) - 
f(W(s_m(r))) \over \mp \Delta x} - f'(W(r \Delta t)) \right| 
< {\epsilon \over 3K}
\label{eq-wdeps} \eeqa 
with probability $1$ assuming $m$ is large enough. 

The function $f'(W(s))$ is continuous with probability $1$, so its 
Riemann sums over $[0, K]$ converge to the corresponding integral as the 
norm of the partition tends to $0$. Thus by (\ref{eq-wdeps}), 
\beqan \lefteqn 
{\left| \int _0^K f'(W(s))~ds -  \sum _{r=1}^{\lfloor K/\Delta t \rfloor } 
{f(W(s_m(r))) - f(W(s_m(r-1))) \over W(s_m(r)) - W(s_m(r-1))} \Delta t 
\right|} \\
& \le &  \sum _{r=1}^{\lfloor K/\Delta t \rfloor } \left| f'(W(r \Delta 
t)) - {f(W(s_m(r)) \mp \Delta x) - f(W(s_m(r))) \over \mp \Delta x} 
\right| \Delta t \\
& + & \left| \int _0^{K_m} f'(W(s))~ds -  
\sum _{r=1}^{\lfloor K/\Delta t \rfloor } f'(W(r \Delta t)) \Delta t 
\right| + \left| \int _{K_m}^K f'(W(s))~ds \right| \\
& < & {\epsilon \over 3K}K + {\epsilon \over 3} + {\epsilon \over 3} 
= \epsilon , 
\eeqan
with probability $1$ if $m$ is large enough. Here $K_m=\lfloor K/\Delta t 
\rfloor \Delta t$. 

Therefore the second sum in (\ref{eq-iteq}) also tends to 
the corresponding integral with probability $1$: 
\[ \lim _{m \to \infty } {1\over 2} \sum _{r=1}^{\lfloor K/\Delta t \rfloor } 
{f(W(s_m(r))) - f(W(s_m(r-1))) \over W(s_m(r)) - W(s_m(r-1))} \Delta t 
= {1\over 2} \int _0^K f'(W(s))~ds . \]

This proves that the defining sum of the It\^o integral in 
(\ref{eq-itod}) converges with probability $1$ as $m \to \infty $, and 
for the limit we have It\^o formula (\ref{eq-itof}). 

Also, by the Stratonovich case of Lemma \ref{l-ditom} and the comments 
made after the lemma, with probability $1$, for all but finitely many $m$, we 
have the following equation for the sum in (\ref{eq-strd}):  
\beqan 
\sum _{r=1}^{\lfloor K/\Delta t \rfloor } {f(W(s_m(r-1))) + 
f(W(s_m(r))) \over 2} (W(s_m(r)) - W(s_m(r-1))) \\
= T_{x=0}^{W(s_m(\lfloor K/\Delta t \rfloor ))}~f(x)\Delta x . 
\eeqan
We saw in (\ref{eq-traco}) that this trapezoidal sum converges to the 
corresponding integral with probability $1$ as $m \to \infty $. Therefore 
the defining sum of the Stratonovich integral in (\ref{eq-strd}) converges 
as well, and for the limit we have formula (\ref{eq-straf}). $\Box $ 

%.............................

\vspace{10pt}

Since the It\^o and Stratonovich formulae are valid for the usual 
definitions of the corresponding stochastic integrals as well, this shows 
that the usual definitions agree with the definitions given in this paper.  

As we mentioned in a special case, the interesting feature of It\^o 
formula (\ref{eq-itof}) is that it contains the non-classical term 
${1\over 2} \int _0^K f'(W(s))~ds$. If $g$ denotes an antiderivative of 
the function $f$, then the It\^o formula can be written as 
\[ g(W(t)) - g(W(0)) = \int _0^t g'(W(s))~dW(s) + {1\over 2} \int 
_0^t g''(W(s))~ds , \]
or formally as the following non-classical chain rule for differentials: 
\[ dg(W(t)) = g'(W(t))~dW(t) + {1\over 2} g''(W(t))~dt . \] 

We mention that other, more complicated versions of It\^o formula can be 
proved by essentially the same method, see \cite{Szabados 1990}. Also, 
as shown there, multiple stochastic integrals can be defined analogously 
as the stochastic integrals defined above. 

%.............................

\vspace{10pt}

{\bf Acknowledgements.} I am obliged to P\'al R\'ev\'esz from whom I learned 
many of the ideas underlying this paper. Also, I am grateful to G\'abor 
Tusn\'ady for his valuable help and to Tam\'as Linder for 
several useful comments. Special thanks are due to 
the Informatics Department of the University of Pisa and particularly 
to Vincenzo Manca for their hospitality while writing the last part of 
the paper.

%-------------------------------------------------------------------------

\end{document}